\documentclass[11pt]{article}
\usepackage{amssymb, amsmath, hyperref}

\def\A{\mathbb{A}}
\def\C{\mathbb{C}}
\def\R{\mathbb{R}}
\def\N{\mathbb{N}}
\def\P{\mathbb{P}}
\def\Q{\mathbb{Q}}
\def\Z{\mathbb{Z}}
\def\J{\mathfrak{J}}
\def\S{\mathcal{S}}

\def\cA {{\cal A}}

\def\cF {{\cal F}}
\def\cM {{\cal M}}
\def\cB {{\cal B}}
\def\cX {{\cal X}}
\def\cR {{\cal R}}
\def\cL {{\cal L}}
\def\cI {{\cal I}}

\def\Res {{\rm Res}}

\newtheorem{defn}{Definition}

\newtheorem{proposition}[defn]{Proposition}
\newtheorem{theorem}[defn]{Theorem}
\newtheorem{corollary}[defn]{Corollary}
\newtheorem{remark}[defn]{Remark}

           {\vspace{3.3mm}
           \noindent{\bf #1}\it}%
           {\vspace{3.3mm}}

\newenvironment{proof}[1]{
  \trivlist \item[\hskip \labelsep{\it #1}]}{\hfill\mbox{$\square$}
  \endtrivlist}

\title{{\bf A parametric representation of totally mixed Nash equilibria}}
\author{Gabriela Jeronimo\thanks{Partially supported
by the following Argentinian research grants: UBACyT X112
(2004-2007) and CONICET PIP 5852/05.} \thanks{Partially supported by
the Argentinian research grant UBACyT X847 (2006-2009).} \quad
Daniel Perrucci$^{{\ast} \dag}$
\\[2mm]
{{\small Departamento de Matem\'atica, Facultad de Ciencias
Exactas y Naturales,}}\\
{{\small Universidad de Buenos Aires, Ciudad Universitaria, 1428
Buenos Aires, Argentina}}\\
{{\small CONICET - Argentina}}\\[3mm] and \\ \ \\ Juan Sabia$^{\ast}$ \\[2mm]
{{\small Departamento de Ciencias Exactas, Ciclo B\'asico Com\'un, and}}\\
{{\small Departamento de Matem\'atica, Facultad de Ciencias
Exactas y Naturales,}}\\
{{\small Universidad de Buenos Aires, Ciudad Universitaria, 1428
Buenos Aires, Argentina}}\\
{{\small CONICET - Argentina}}}
\date{}

\begin{document}

\maketitle

\begin{abstract}

We compute a parametric description of the totally mixed Nash
equilibria of a generic game in normal form with pre-fixed
structure. Using this representation, we show 
conditions under which a game has the maximum possible number of
this kind of equilibria. Then, we present a symbolic procedure
that allows us to describe and estimate the number of isolated
totally mixed Nash equilibria of an arbitrary game. Under certain
assumptions, the algorithm computes the exact number of these
equilibria.
\end{abstract}

\section{Introduction}

Noncooperative game theory is used to model and analyze strategic
interaction situations. Among its most outstanding applications,
we can mention the fundamental role this theory has played in
economics (see, for example, the classical reference book
\cite{vonNeumannMorg}). Moreover, game theory has also been
applied to politics, sociology and psychology, and to biology and
evolution as well.

One of the main concepts in this theory is a \emph{Nash
equilibrium}, which consists in a situation in which no player can
increase his payoff by unilaterally changing his strategy. Since
within this theory the players cannot communicate in order to
decide a simultaneous change of strategies, in a Nash equilibrium
the game stabilizes. In \cite{Nash}, it is proved that any
noncooperative game in normal form has at least one Nash
equilibrium. However, the proof is not constructive and does not
give any information about the existence of more than one Nash
equilibrium. The question posed is how to compute algorithmically
Nash equilibria and to determine the number of them in a given
game.

Nash equilibria of noncooperative games in normal form can be
regarded as real solutions to systems of polynomial equations and
inequalities (see, for instance,  \cite[Chapter 6]{Stu02}). In the
case of two players, each equilibrium is the solution of a linear
system of equations, and therefore, equilibria may be found
exactly by using simplex type algorithms (see, for instance,
\cite{LemkeH}); however, there is no polynomial time algorithm
solving the problem (see \cite{vSt02}). In the general case of a
game with more than two players, the polynomials appearing are
multilinear. To solve the problem of finding \emph{one}
equilibrium, some numerical methods have been applied successfully
(for example, some methods derived from Scarf's algorithm,
\cite{Sca67}). Nevertheless, sometimes it is not sufficient to
compute only one equilibrium because, depending on the problem to
be solved, not all the equilibria of a game are equally
interesting and the methods developed to compute only one
equilibrium do not allow us to decide beforehand whether it
fulfills some additional properties or to compare different
equilibria.

A comparative study of different known methods for the computation
of all the Nash equilibria of a game may be found in \cite{Datta}.
In \cite{HP05}, a new algorithm solving this problem for generic
games by means of homotopy methods is presented, but no complexity
bounds are shown. In addition, the application of symbolic
algorithms solving systems of equations and inequalities over the
real numbers is being studied in this context, motivated by the
characterization of the set of all the Nash equilibria of a game
as a semi-algebraic set (an example of this fact is the
application of quantifier elimination algorithms over the real
numbers to compute approximated equilibria in \cite{LM04}; see
also the survey \cite{McKMcL96}). However, up to now, no
significative result had been obtained concerning the adaptation
of these algorithms in order to profit from the particular
properties of the algebraic systems arising in game theory.

In this paper, we study \emph{totally mixed Nash equilibria}, that
is to say, Nash equilibria in which every player allocates a
positive probability to each of his available strategies. Note
that a procedure to compute these equilibria can be used as a
subroutine to compute all Nash equilibria of the game by recursing
over all possible subsets of used strategies. We present a
\emph{symbolic} method to find a parametric description of the set
of totally mixed Nash equilibria of a generic game with a
pre-fixed structure. This method is based on a symbolic procedure
for the computation of \emph{multihomogeneous} resultants with
complexity polynomial in the degree and the number of variables of
the resultant (\cite{JeSa}). Using the description previously
obtained, we show conditions (given by polynomial inequalities on the
payoff values) under which a game with
the given pre-fixed structure will have the maximum number of such
equilibria. The next step is to give a similar parametric
description for \emph{particular} games. First, we solve this
problem under some genericity assumptions implying, in particular,
that the number of totally mixed Nash equilibria of the game is
finite, and we show how to compute this number. Then, we consider
the general case, in which we give a parametric description of a
finite set of points including all the isolated (in the complex
space) totally mixed Nash equilibria of the game, which enables us
to bound the number of these isolated equilibria.
 All our algorithms
have a complexity polynomial in the number of pure strategies of
each player, and the number of totally mixed Nash equilibria of a
generic game with the given structure.

This paper is organized in the following way:

In Section \ref{Preliminaries}, we introduce some basic notions on
game theory, the description of the algorithmic model we are going
to use and a mathematical formulation  for the totally mixed Nash
equilibria of a game in normal form as the set of solutions of a
system of multihomogeneous polynomial equations. In Section
\ref{sec generic}, the algorithmic description of the totally
mixed Nash equilibria of a generic game is given, and it is used
to find conditions under which a game has the maximum number of
these equilibria. In Section \ref{sec particular} we deal with the
totally mixed Nash equilibria of particular games. Finally, the
last section is devoted to proving some results about the
algorithmic computation of multihomogeneous resultants and upper
bounds for their degrees which are used throughout the paper.

\section{Preliminaries}\label{Preliminaries}

\subsection{Definitions and Notation}

\subsubsection{Basic facts}

Throughout this paper $\mathbb{Q}$ denotes the field of rational
numbers, $\mathbb{N}$ denotes the set of positive integers and
$\mathbb{N}_0:= \mathbb{N} \cup \{ 0 \}$.

If $K$ is a field, we denote an algebraic closure of $K$ by
$\overline K$. The ring of polynomials in the variables
$x_1,\dots, x_n$ with coefficients in $K$ is denoted by
$K[x_1,\dots, x_n]$. For a polynomial $f \in K[x_1,\dots, x_n]$ we
write $\deg f$ to refer to the total degree of $f$ and
$\deg_{x_i}f$ to refer to the degree of $f$ in the variable $x_i$.

For $n\in \mathbb{N}$ and an algebraically closed field $k$, we
denote by $\mathbb{A}^n(k)$ and $\mathbb{P}^n(k)$ (or simply by
$\mathbb{A}^n$ or $\mathbb{P}^n$ if the base field is clear from
the context) the $n$-dimensional affine space and projective space
over $k$ respectively, equipped with their Zariski topologies. 

We adopt the usual notions of dimension and degree of an
algebraic variety $V$, which will be denoted by $\dim V$ and
$\deg V$ respectively. See for instance \cite{Shafarevich} and
\cite{Heintz83} for  the definitions of these notions.

\subsubsection{Game theory}

In this section we present some basic game theory concepts. For a
more detailed account on the subject we refer the reader to any
standard game theory text such as \cite{OR94}.

We consider non-cooperative games in \emph{normal form}; that is
to say, games in which there is only one time step at which all
the players move simultaneously without communicating among
themselves. We will assume that there are $r$ players in the game
having $n_1 +1, \dots, n_r +1$ distinct available pure strategies
respectively ($n_1, \dots, n_r \in \N$).

For $i=1, \dots, r$:
\begin{itemize}

\item $c^{(i)}:= ( c^{(i)}_{j_1 \dots j_r})_{0\le j_k \le n_k}$ is
the given \emph{payoff matrix} of  player $i$, where $c^{(i)}_{j_1
\dots j_r}$ is the \emph{payoff} to player $i$ if, for every $1\le
k \le r$, player $k$ chooses the strategy $j_k$.

\item $X_i:=(x_{i0}, x_{i1}, \dots, x_{in_i})$ is a vector representing a
\emph{mixed strategy} of the $i$th player, which is a probability
distribution
 on his set of pure strategies
 (that is to say, $x_{ij}$ is the probability that the $i$th player
  allocates to his $j$th pure strategy).

\end{itemize}
With these notations, for every $1\le i \le r$, the payoff to
player $i$ if the mixed strategies $X_1,\dots, X_r$ are played is:
$$\pi_i(X_1,\dots, X_r):= \sum_{0\le j_1\le n_1}\dots \sum_{0\le
j_r\le n_r} c^{(i)}_{j_1 \dots j_r} x_{1 j_1} \dots x_{r j_r}.$$

A \emph{Nash equilibrium} is a vector of mixed strategies such
that no player can increase his payoff by changing unilaterally to
another strategy while the other players keep their strategies
fixed; that is to say, a vector of mixed strategies $X_1,\dots,
X_r$ satisfying
$$
\pi_i(X_1,\dots, X_{i-1},X_i, X_{i+1},\dots, X_r) \ge \pi_i(X_1,
\dots, X_{i-1}, X'_i, X_{i+1}, \dots, X_r) \quad \forall\, 1\le i
\le r
$$
for every mixed strategy $ X'_i$. A \emph{totally mixed Nash
equilibrium} is a Nash equilibrium in which each pure strategy is
assigned a positive probability, that is, one that satisfies
$x_{ij}>0$ for every $1\le i \le r, \ 0\le j\le n_i$.

\subsection{Data structure, algorithms and complexity}
\label{data structures}

The algorithms we consider in this paper are described by arithmetic
networks over the base field $\mathbb{Q}$ (see \cite{vzg86}). An
arithmetic network is represented by means of a directed acyclic
graph. The external nodes of the graph correspond to the input and
output of the algorithm. Each of the internal nodes of the graph is
associated with either an arithmetic operation in $\mathbb{Q}$ or a
comparison between two elements in $\mathbb{Q}$ followed by a
selection of another node. These are the only operations allowed in
our algorithms.

We assume that the cost of each operation in the algorithm is 1
and so, we define the {\em complexity} of the algorithm as the
number of internal nodes of its associated graph.

The objects our algorithms deal with are polynomials with
coefficients in $\mathbb{Q}$. We represent each of them by means of
one of the following data structures:
\begin{itemize}
\item {\em Dense form}, that is, as the array of all its
coefficients (including zeroes) in a pre-fixed order of all
monomials of degree at most $d$, where $d$ is an upper bound for
the degree of the polynomial. The size of this representation
equals the number of monomials of degree at most $d$. \item {\em
Sparse encoding,} that is, as an array of the coefficients
corresponding to monomials in a fixed set, provided that we know
in advance that the coefficient of any other monomial of the
polynomial must be zero. The size in this case is the cardinal of
the fixed set of monomials. \item {\em Straight-line programs,}
which are arithmetic circuits (i.e. networks without branches).
Roughly speaking, a straight-line program (or slp, for short) over
$\mathbb{Q}$ encoding a polynomial $f\in \mathbb{Q}[x_1,\dots,
x_n]$ is a program which enables us to evaluate the polynomial $f$
at any given point in $\mathbb{Q}^n$. Each of the instructions in
this program is an addition, subtraction or multiplication between
two pre-calculated elements in $\mathbb{Q}[x_1,\dots, x_n]$, or an
addition or multiplication by a scalar. The number of instructions
in the program is called the {\em length} of the slp. For a
precise definition of straight-line programs we refer to
\cite[Definition 4.2]{Burgisser}  (see also \cite{HS82}).
\end{itemize}

\subsection{Polynomial equations for totally mixed Nash equilibria}\label{TMNEeq}

Now, we show that the totally mixed Nash equilibria of an
$r$-person game in normal form can be regarded as real solutions
to a polynomial equation system (see \cite[Sec.~6.3]{Stu02}).

Let us observe that, if the payoff matrix of player $i$ is
$c^{(i)}$ ($1 \le i \le r)$, given mixed strategies $X_1,\dots,
X_r$,  the payoff to player $i$ is
\begin{eqnarray*}
\pi_i(X_1,\dots, X_r)&=& \sum_{j_1, \dots, j_r} c^{(i)}_{j_1
\dots j_r} x_{1 j_1} \dots x_{r j_r}\\
&=&\sum_{0\le j_i\le n_i} x_{ij_i}\sum_{j_1,\dots,
j_{i-1},j_{i+1},\dots, j_r} c^{(i)}_{j_1 \dots j_r} x_{1 j_1}
\dots x_{i-1 j_{i-1}}x_{i+1 j_{i+1}}\dots x_{r
j_r}\\
&=& \sum_{0\le j_i\le n_i} x_{ij_i}\, \pi_i(X_1,\dots, X_{i-1},
e^{(i)}_{j_i}, X_{i+1}, \dots, X_r)
\end{eqnarray*}
where $e_j^{(i)}$ is the $(j+1)$th vector of the canonical basis of
$\mathbb{R}^{n_i + 1}$.

Now, if $(X_1,\dots, X_r)$ is a totally mixed Nash equilibrium, we have that $x_{ij}
>0$ for every $i,j$. Then, in order that
$\pi_i(X_1,\dots, X_{i-1}, \cdot\, , X_{i+1}, \dots, X_n)$ attains
a maximum at $X_i = (x_{i0}, \dots, x_{in_i})$ (among all the
different probability distributions), a necessary and sufficient
condition is that
\begin{equation}\label{indiference}
\pi_i(X_1,\dots, X_{i-1}, e_{j_i}^{(i)}, X_{i+1}, \dots, X_r)  =
\pi_i(X_1,\dots, X_{i-1}, e_{0}^{(i)}, X_{i+1}, \dots, X_r) \quad
\forall\, 1\le j_i \le n_i.
\end{equation}

In fact, assuming that this condition does not hold for some
$1\le i \le r$, let $j_{\rm max}$ and $j_{\rm min}$ denote the
indices of pure strategies leading to maximum and minimum payoffs
to player $i$ respectively. Then, taking $\varepsilon
>0$ with $\varepsilon <\min\{x_{i j_{\rm min}}, 1-x_{i j_{\rm max}}\}$
and changing to the (totally) mixed strategy $X'_i := X_i +
\varepsilon (e^{(i)}_{j_{{\rm max}}} - e^{(i)}_{j_{{\rm min}}})$
the $i$th player can increase his payoff by $$\varepsilon
(\pi_i(X_1,\dots, X_{i-1}, e^{(i)}_{j_{\rm max}}, X_{i+1}, \dots,
X_r) - \pi_i(X_1,\dots, X_{i-1}, e^{(i)}_{j_{\rm min}}, X_{i+1},
\dots, X_r))>0,$$ contradicting the fact that the vector is a
totally mixed Nash equilibrium for the game.

Conversely, for a vector of totally mixed strategies $(X_1,\dots,
X_r)$ such that all the identities (\ref{indiference}) hold for
every $1\le i \le r$, we have that $\pi_i(X_1,\dots,  X'_i, \dots,
X_r) = \pi_i(X_1,\dots, X_i, \dots, X_r)$ for every mixed
strategy $X'_i$ and so, player $i$ $(1\le i \le r)$ cannot
increase his payoff by unilaterally changing his strategy, which
implies that $(X_1,\dots, X_r)$ is a Nash equilibrium for the game.

We conclude that the totally mixed Nash equilibria of the game are those vectors
$(X_1,\dots, X_r)$ with $X_i := (x_{i0}, \dots, x_{in_i})$ for
every $1\le i \le r$ satisfying:
\begin{itemize}
\item[(a)] $x_{ij} > 0$ \ for $i=1,\dots, r$ and $j=0,\dots, n_i$,

\item[(b)] $\displaystyle\sum_{0\le j \le n_i} x_{ij} = 1$ \ for $i=1,\dots, r$,

\item[(c)] $\displaystyle\sum_{J_{-i}}
 \ a^{(ik)}_{J_{-i}} \ x_{1j_1} \dots
x_{{i-1}j_{i-1}}x_{{i+1}j_{i+1}} \dots x_{rj_r} = 0$ \  for $ i
=1,\dots, r$ and $k=1,\dots, n_i$, where the sum runs over all
$J_{-i} := j_1\dots j_{i-1} j_{i+1} \dots j_r$ with $0\le
j_t \le n_t$ for every $t\ne i$ and  $a^{(ik)}_{J_{-i}} :=
c^{(i)}_{j_1\dots j_{i-1} k j_{i+1} \dots j_r} -
c^{(i)}_{j_1\dots j_{i-1} 0 j_{i+1} \dots j_r}$.
\end{itemize}

Let us observe that (c) is a system of $n:= n_1 + \dots + n_r$
multihomogeneous polynomial equations in the $r$ groups of
variables $X_1,\dots, X_r$ with $n_1+1, \dots, n_r+1$ variables
respectively (with degrees $1$ or $0$ with respect to each group)
and, therefore, it defines a (possibly empty) projective variety
in $\P^{n_1}(\mathbb{C}) \times \dots \times
\P^{n_r}(\mathbb{C})$. The complex solutions to the polynomial
equation system (c) will be called \emph{quasi-equilibria} of the
game (see \cite{Datta}), and those solutions not lying in any of
the infinite hyperplanes $\{x_{i0} = 0\}$ $(1\le i \le r)$ will be
called \emph{affine quasi-equilibria} of the game.

Every quasi-equilibrium $\xi:= (\xi_1, \dots, \xi_r)\in
\P^{n_1}(\mathbb{C}) \times \dots \times \P^{n_r}(\mathbb{C})$
determines  at most one totally mixed Nash equilibrium of the
game: for every $1\le i \le r$, let $s_{\xi_i}:= \sum_{0\le j \le
n_i} \xi_{ij}$ be the sum of the coordinates of $\xi_i$. If
$s_{\xi_i} \ne 0$ for every $1\le i \le r$, the unique associated
representation of $\xi$ whose coordinates satisfy condition (b) is
$(\xi_1/s_{\xi_1}, \dots, \xi_r/s_{\xi_r})$, and it will be a
totally mixed Nash equilibrium if and only if all its coordinates
are positive real numbers.

\subsection{On the number of solutions to a multihomogeneous
system}\label{number of solutions}

Let $r\in \mathbb{N}$ be a positive integer. Fix positive integers
$n_1, \dots, n_r$ and consider $r$ groups of variables $X_j:=
(x_{j0},\dots, x_{j\, n_j})$, $j=1,\dots, r$. We say that the
polynomial $F \in K[X_1,\dots, X_r]$ is {\em multihomogeneous} of
{\em multi-degree} $v:=(v_1,\dots, v_r)$, where $v_1,\dots, v_r$
are non-negative integers, if $F$ is homogeneous of degree $v_j$
in the group of variables $X_j$ for every $1\le j \le r$.

Set $n:= \sum_{j=1}^r n_j$. The classical Multihomogeneous
B\'ezout Theorem, which follows from the intersection theory for
divisors  (see for instance \cite[Chapter 4]{Shafarevich}), states
that the set of common zeroes (over an algebraically closed field)
in the projective variety $\mathbb{P}^{n_1} \times \dots \times
\mathbb{P}^{n_r}$ of $n$ generic multihomogeneous polynomials
$F_1,\dots, F_n$ with multi-degrees $\nu_i:=(\nu_{i1}, \dots,
\nu_{ir})$ for $i=1,\dots, n$ is a zero-dimensional variety with
\begin{equation}\label{bezoutmultihom}
{\rm Bez}_{n_1\dots, n_r}(\nu_1;\dots; \nu_n) :=
\sum_{(j_1,\dots, j_n)\in \J} \Big(\prod_{i=1}^n \nu_{i
j_{i}}\Big)\end{equation} points, where $\J:= \{ (j_1,\dots, j_n)
\ /\ \#\{ k : j_k = i \} = n_i \ \forall \, 1\le i \le r\}$. For
an alternative proof of this result using deformation techniques,
we refer the reader to  \cite{MSW95}, \cite{McLennan} and
\cite{HJSS05}. Note that this can be seen as a particular case of
Bernstein's theorem on the number of common roots of sparse
systems \cite{Bernstein}.

The quantity ${\rm Bez}_{n_1\dots, n_r}(\nu_1;\dots; \nu_n)$ is
called the \emph{B\'ezout number} of the generic multihomogeneous
polynomial system. If $k_1,\dots, k_t$ are positive integers with
$\sum_{i=1}^t k_i = n$, we will use the notation ${\rm
Bez}_{n_1,\dots, n_r}(\nu_1, k_1; \dots ; \nu_t, k_t)$ for the B\'
ezout number of a multihomogeneous system with $k_i$ polynomials
of multi-degree $\nu_i$ for every $1\le i \le t$.

The equations arising in our particular setting for the
computation of totally mixed Nash equilibria of a game are
multilinear (see Section \ref{TMNEeq}). Moreover, for a game with
$r$ players with $n_1+1, \dots, n_r+1$ pure strategies
respectively, we have a system of $n= \sum_{j=1}^r n_j$ polynomial
equations consisting of exactly $n_i$ polynomials of multi-degrees
equal to $d_i := (1,\dots, 1,0,1,\dots, 1) \in (\N_{0})^r$ (where
the $0$ lies in the $i$th coordinate of $d_i$) for every $1\le i
\le r$. Then, according to the previous formula, a system coming
from a generic game with the given structure will have
\begin{equation}\label{deltacalculo}
\delta := \sum_{(j_{11},\dots,j_{1n_1},\dots, j_{r1},\dots,
j_{rn_r})\in \J} \Big(\prod_{i=1}^r \prod_{k=1}^{n_i} d_{i
j_{ik}}\Big) \end{equation} solutions in $\mathbb{P}^{n_1} \times
\cdots \times \mathbb{P}^{n_r}$, where $\J:= \{ (j_{11},\dots, j_{1
n_1},\dots, j_{r1},\dots, j_{rn_r}) \ /\ \#\{ j_{hk} : j_{hk}= i
\} = n_i \ \forall \, 1\le i \le r\}$. In fact, for a ``generic''
game, this will be the number of totally mixed Nash equilibria
(see \cite{McKelvey}).

We are going to deal with the case in which $\delta >0$. This
inequality can be determined by considering the set of exponents
appearing with non-zero coefficients in each of the polynomials in
the system (see \cite[Chapter IV, Proposition 2.3]{Oka}) and in
our particular case, it is equivalent to the fact that $n_j \le
\sum_{1\le k \le r,\ k\ne j} n_k = n- n_j$ for every $1\le j \le
r$. From now on, we will assume that these inequalities hold.

Taking into account that $d_{ii} = 0$ for every $1\le i \le r$,
the only $n$-tuples contributing to the sum (\ref{deltacalculo})
are those where $j_{i k} \ne i $ for every $1\le k \le n_i$, and
for each of them, the corresponding term of the sum equals $1$.
Therefore, $\delta$ equals the cardinality of the set
\begin{equation}\label{J0}
\J_0 = \{ (j_{11},\dots, j_{rn_r}) \, /\, j_{ik} \ne i \,
\forall\, 1\le k \le n_i \hbox{ and } \#\{ j_{hk} \, /\, j_{hk}=
i \} = n_i \ \forall \, 1\le i \le r\}.
\end{equation}

\subsection{Geometric resolutions}

A way of representing zero-dimensional affine varieties which is
widely used in computer algebra nowadays is  a \emph{geometric
resolution}. This notion was first introduced in the works of
Kronecker and K{\"o}nig in the last years of the XIXth century
\cite{Kronecker}. Roughly speaking, it consists in a parametric
description of the variety in which the parameter values range
over the set of roots of a univariate polynomial. Now, we give the
precise definition we are going to use.

Let $V = \{ \xi^{(1)}, \dots,\xi^{(\delta)} \} \subset \A^n$ be a
zero-dimensional variety defined by polynomials in $K[x_1, \dots,
x_n]$. Given a \emph{separating} linear form $\ell = u_1 x_1 +
\dots + u_n x_n\in K[x_1, \dots, x_n]$ for $V$ (that is, a linear
form $\ell$ such that $\ell (\xi^{(i )})\ne \ell(\xi^{(k)})$ if
$i \ne k$), the following polynomials completely characterize the
variety $V$:
\begin{itemize}
\item the \emph{minimal polynomial} $p:= \prod_{1 \le i \le
\delta} (T - \ell(\xi^{(i)}))\in K[T]$ of $\ell$ over the variety
$V$ (where $T$ is a new variable), \item polynomials $w_1,\dots,
w_n \in K[T]$ with $\deg w_j < \delta$ for every $1\le j \le n$
satisfying
$$ V = \big\{ \big(\frac{w_1}{p'}(\eta) ,\dots, \frac{w_n}{p'} (\eta) \big) \in \overline{K}^n
/ \, \eta \in \overline{K} ,\ p(\eta) = 0 \big\}.$$
\end{itemize}
The family of univariate polynomials $p, w_1,\dots, w_n\in K[T]$
is called the \emph{geometric resolution} of $V$ (associated with
the linear form $\ell$).

In our particular setting of totally mixed Nash equilibria
computation, we will not only deal with zero-dimensional
varieties in an affine space, but we will also consider
zero-dimensional subvarieties of $\P^{n_1} \times \dots \times
\P^{n_r}$.

Write $\xi:=(\xi_1,\dots, \xi_r)\in \P^{n_1} \times \dots \times
\P^{n_r}$ with $\xi_i:=(\xi_{i0}: \dots : \xi_{in_i})$ $(1\le i
\le r)$ to denote a point in the projective variety  $\P^{n_1}
\times \dots \times \P^{n_r}$. Assume that $V\subset \P^{n_1}
\times \dots \times \P^{n_r}$ is a zero-dimensional  variety
defined by multihomogeneous polynomials in $K[X_1,\dots, X_r]$
such that $\xi_{i0} \ne 0$ $(1\le i \le r)$ holds for every point
$\xi \in V$. Then, we may associate with $V$ the following
zero-dimensional variety in $\A^n$, where $n:= n_1+\cdots +n_r$:
$$\{(\xi'_1,\dots, \xi'_r)\in \A^n \,/\, \xi'_i =
(\xi_{i1}/\xi_{i0},\dots,\xi_{in_i}/\xi_{i0}) \ \forall\, 1\le i
\le r, \ \xi \in V\}.$$ A geometric resolution $p, w_{11},\dots,
w_{1n_1},\dots, w_{r1},\dots, w_{rn_r}\in K[T]$ of this
zero-dimensional variety will also be called a \emph{geometric
resolution} of $V\subset \P^{n_1} \times \dots \times \P^{n_r}$.
In this case, the geometric resolution of $V$ provides the
following description of the variety:
$$V = \big\{ \big((p'(\eta):w_{11}(\eta):\dots : w_{1n_1}(\eta)), \dots,
(p'(\eta):w_{r1}(\eta):\dots : w_{rn_r}(\eta))\big)  /\, \eta\in
\overline K,\ p(\eta)= 0 \big\}.$$

\section{On the totally mixed Nash equilibria of a generic
game}\label{sec generic}

This section is devoted to the study of totally mixed Nash equilibria of generic games.
In order to do this, we will begin by treating the payoff values as parameters and  computing
a geometric resolution of the set of quasi-equilibria of the associated generic game. This geometric resolution will provide a rational formula of the parameters which, for generic values of them, represents the totally mixed Nash equilibria of the game.

Then, we will obtain a finite number of generic conditions ensuring that a specific game has the maximum possible number of totally mixed Nash equilibria. Under these conditions, the Nash equilibria of the game can be described by the geometric resolution obtained by substituting the given payoffs for the parameters in our generic formulas.

\subsection{The set of quasi-equilibria of a generic game}
\label{casogen}

Here we are going to compute a geometric resolution of the set of quasi-equilibria
of a generic game with $r$ players with $n_1+1, \dots, n_r +1$ pure strategies respectively,
where $n_i \in \N$ for every $1\le i \le r$.

For $1 \le i \le r, \, 1 \le k \le n_i$, let $A^{(ik)} := (
A^{(ik)}_{j_1 \dots j_{i-1} j_{i+1} \dots j_r})_{0\le j_t \le
n_t}$ be a set of new indeterminates and
\begin{equation}\label{genpoly}
F^{(i)}_{k} := \sum_{J_{-i}}
 \ A^{(i k)}_{J_{-i}} \ x_{1j_1} \dots
x_{{i-1}j_{i-1}}x_{{i+1}j_{i+1}} \dots x_{rj_r} ,
\end{equation}
where the sum runs over all $J_{-i} := j_1\dots j_{i-1}
j_{i+1} \dots j_r$ with $0\le j_t \le n_t$ for every $t\ne i$;
that is, $F^{(i)}_k$ is a generic multihomogeneous polynomial of
multi-degree $d_i:=(1,\dots, 1, 0, 1,\dots, 1)\in (\N_{0})^r$
(where $0$ is in the $i$-th coordinate).

Note that, if $(C^{(i)}_{j_1\dots  j_r})_{1\le i \le r,\, 0\le j_t
\le n_t}$ are new variables that are regarded as parameters
representing the generic payoffs of the game, from the polynomials
introduced in (\ref{genpoly}), we can obtain the polynomials
defining the quasi-equilibria of the generic game
 by substituting
\begin{equation}\label{subsgame}
A^{(i k)}_{J_{-i}} := C^{(i)}_{j_1\dots j_{i-1} k j_{i+1} \dots
j_r} - C^{(i)}_{j_1\dots j_{i-1} 0 j_{i+1} \dots j_r}
\end{equation}
for every $1\le i \le r$, $1\le k\le n_i$ and each $J_{-i} =
j_1\dots j_{i-1} j_{i+1} \dots j_r$ with $0\le j_t \le n_t$
for every $t\ne i$.
 Thus,
we are going to compute a geometric resolution of the solution
set of the generic polynomial system $F^{(i)}_{k} = 0$, $1\le i
\le r,\, 1\le k \le n_i,$ and then, we will make the
substitution (\ref{subsgame}) in the result to get the desired
geometric resolution of the quasi-equilibria of the generic game.

Our main result is the following:

\begin{theorem}\label{eq-quasieq} Let notations be as in Section \ref{Preliminaries}.
There is an algorithm which computes a geometric resolution of the
set of quasi-equilibria of a generic game with $r$ players having
$n_1 +1, \dots, n_r+1$ pure strategies respectively, within
complexity
 $O( D^2 (D + n_1\dots n_r\,\delta \log(D)
r^2 n^4 (n^3 + r N)))$, where
$$\begin{array}{rcl}
D &:=& \sum\limits_{0\le i \le r} n_i{\rm Bez}_{n_1,\dots,
n_r}(d_0,n_0; d_{1}, n_{1}; \dots; d_{i},n_i-1; \dots; d_{r}, n_r) \quad (n_0 \! := 1, d_0\! :=(1,\dots,1)), \\[4mm]
\delta &:=& {\rm Bez}_{n_1,\dots, n_r}
(d_{1}, n_{1}; \dots; d_{r}, n_r),\\[3mm]
n&:= & \sum\limits_{1\le i \le r} n_i,\\[3mm]
 N &:=& n + 1 + \sum\limits_{1\le i \le r} n_i(n_1 + 1)\dots(n_{i-1} + 1)(n_{i+1} +
1)\dots(n_r + 1). \end{array}$$ The algorithm obtains polynomials
$P(T),\, W_{ij}(T) \in \Q[C^{(i)}_{j_1\dots  j_r}][T]$ with
$\deg_TP = \delta$, $\deg_T W_{ij}< \delta $ and degrees bounded
by $D$ in the parameters $C^{(i)}_{j_1\dots j_r}$, which are
encoded by straight-line programs of length $O( D^2 (D + n_1\dots
n_r\,\delta \log(D) r^2 n^4 (n^3 + r N)))$.
\end{theorem}

\begin{proof}{Proof.}
The geometric resolution of the system $F^{(i)}_{k} = 0$, $1\le i
\le r, 1\le k \le n_i,$ will be obtained by means of
multihomogeneous resultant computations. More precisely, we are
going to obtain the geometric resolution  associated to a linear
form $\ell$ which separates the points of the generic system by
computing the resultant of our system and a generic linear form.
The minimal polynomial of $\ell$ and the polynomials $W_{ij}$ will
be obtained by specializing the resultant and some of its partial
derivatives in the coefficients of $\ell$.

We introduce a set of new indeterminates $A^{(0)}:= (A^{(0)}_0,
A^{(0)}_{ij} : 1\le i \le r,\, 1\le j \le n_i)$ which stand for
the coefficients of a generic affine linear form $A^{(0)}_0
+\sum_{1\le i \le r \atop 1\le j \le n_i} A^{(0)}_{ij} \, x_{ij}$
in the $n= n_1 +\cdots +n_r$ variables $x_{ij}$ $(1\le i \le r,
1\le j \le n_i)$, and we consider the multilinear polynomial
$$F_0 := A^{(0)}_{0} x_{10} \dots  x_{r0} +
\sum_{\tiny\begin{array}{c} 1\le i\le r\\
1\le j \le n_i\end{array}}  A^{(0)}_{ij} \, x_{10} \dots x_{i-1
\, 0}x_{i j}x_{i+1\, 0} \dots x_{r0},$$ which is obtained by
homogenizing  the generic affine linear form  with respect to
each group of variables $X_1,\dots, X_r$.

The first step of the algorithm consists in the computation of the
polynomial resultant $\mathcal{R} = {\rm Res}(F_0, F_1^{(1)},
\dots, F_{n_1}^{(1)}, \dots, F_1^{(r)}, \dots, F_{n_r}^{(r)})$.
According to \cite[Theorem 1]{Min}, this resultant coincides with
the specialization of the resultant of a system of
multihomogeneous polynomials of respective multi-degrees $d_0 =
(1, \dots, 1), d_{j}^{(1)} = (0, 1, \dots, 1)$ for $1 \le j \le
n_1,
 \dots, d_{j}^{(r)} = (1, \dots,
1, 0)$ for $1 \le j \le n_r$, in which all the coefficients of the
polynomial of multi-degree $d_0$ corresponding to monomials not
appearing in $F_0$ are substituted for $0$. Thus, in order to
obtain the polynomial $\mathcal{R}$ we apply an adapted version of
the algorithm underlying the proof of \cite[Theorem 5]{JeSa} (see
Subsection \ref{appendix}). With our previous notation, the
complexity of this algorithm is of order $O( D^2 (D + n_1\dots
n_r\delta \log(D) r^2 n^4 (n^3 + r N)))$ and it computes a
straight-line program for $\mathcal{R}$ whose length is of the
same order.

Then, a standard procedure enables us to compute a geometric
resolution of the zero-dimensional variety defined by the system
$F^{(i)}_k = 0$ $(1\le i \le r, \, 1\le k \le n_i)$ from the
polynomial $\mathcal{R}$: the polynomials giving the
parametrization of the points in the variety are obtained from the
partial derivatives
$$\mathcal{R}_0:= \frac{\partial \mathcal{R}}{\partial A^{(0)}_0} \quad \hbox{ and
} \quad \mathcal{R}_{ij}:= \frac{\partial \mathcal{R}}{\partial
A^{(0)}_{ij}}\  (1 \le i \le r,\, 1 \le j \le n_i),$$  which we
compute from the straight-line program representing $\mathcal{R}$
within the same complexity order as in the first step. These
partial derivatives are encoded by straight-line programs of
length of the same order as for the straight-line program encoding
$\mathcal{R}$ (see \cite{Burgisser}).  Note that $\deg_{A^{(0)}_0}
\mathcal{R}_{ij}<\delta$ for every $1\le i \le r, 1 \le j \le
n_i$.

Let $L:= \sum_{1\le i \le r\atop 1 \le j \le n_i} L_{ij} x_{ij}$
be a generic linear form in the variables $x_{ij}$, where $L_{ij}$
are new variables, and let $T$ be another new variable. Consider
the polynomial $P_L \in \mathbb{Q}[A^{(ik)}, L_{ij}][T]$ obtained
by specializing
\begin{equation}\label{subst}
A_0^{(0)} \mapsto T, \ A_{ij}^{(0)} \mapsto  - L_{ij}\ (1\le i \le
r, 1 \le j \le n_i)
\end{equation}
 in $\mathcal{R}$. Since the resultant
$\mathcal{R}$ is in the ideal $(F_0, F^{(i)}_{k}: 1\le k \le r,
1\le i\le n_k)$, then substituting $L$ for $T$ in $P_L$, we obtain
a polynomial $\mathcal{P}$ which is in the ideal $(F^{(i)}_{k}:
1\le k \le r, 1\le i\le n_k)$. Taking into account that
$\deg_T(P_L) = \delta$, which is the number of solutions of the
system $F_k^{(i)}=0$, it follows that $P_L$ is a multiple by a
nonzero factor in $\mathbb{Q}[A^{(ik)}]$ of the minimal polynomial
of $L$. On the other hand, taking the partial derivative of
$\mathcal{P}$ with respect to the variable $L_{ij}$, we get
$$ \frac{\partial \mathcal{P}}{\partial L_{ij}} = -\frac{\partial \mathcal{R}}{\partial A^{(0)}_{ij}}(L, -L_{ij}, A^{(ik)})+
\frac{\partial \mathcal{R}}{\partial A^{(0)}_{0}}(L, -L_{ij},
A^{(ik)}) x_{ij} \in (F^{(i)}_{k}: 1\le k \le r, 1\le i\le n_k).
$$
We conclude that making the substitution (\ref{subst}) in
$\mathcal{R}_0$ and $\mathcal{R}_{ij}$ $(1\le i \le r, 1 \le j \le
n_i)$ we obtain polynomials which complete the geometric resolution
of the variety defined by $F^{(i)}_k = 0$ with respect to the
generic linear form $L:= \sum_{1\le i \le r\atop 1 \le j \le n_i}
L_{ij} x_{ij}$.

Finally, in order to obtain a geometric resolution for the zero
set of the generic system, we choose a separating linear form and
we substitute its coefficients for the parameters $L_{ij}$. As the
multihomogeneous system $F^{(i)}_k = 0$ is generic, it has no
zeroes in the hyperplanes $x_{i0} =0$ $(1 \le i \le r)$ and we can
consider its zeroes as affine points by setting $x_{i0} =1$ $(1
\le i \le r)$. Let us show that the linear form $\ell: =
\sum_{1\le i \le r} x_{i1}$ separates these affine points. To this
end, it is enough to show the existence of a choice of
coefficients such that the induced polynomial system has the
maximum number of affine roots and that $\ell$ separates those
roots.  To see this, choose coefficient vectors for the
polynomials  $F^{(i)}_k$ so as to obtain a specific system
$f^{(i)}_k$ with the maximum number of affine solutions, and take
a linear form $l\in \Q[x_{ij}; 1\le i \le r, 1\le j \le n_i]$
separating the common affine roots of the polynomials $f^{(i)}_k$.
Now, making a linear change of variables in each group $X_i$ $(1
\le i \le r)$, the linear form $l$ maps to $\ell$ and the specific
system considered leads to a polynomial system of the same
structure in the new variables, which is the particular system we
were looking for.

Hence, the algorithm proceeds to specialize
$$A_0^{(0)} \mapsto T, \ A_{i1}^{(0)} \mapsto -1 \ (1\le i \le r), \ A_{ij}^{(0)} \mapsto
0 \ (1\le i \le r, 2 \le j \le n_i)$$ in the polynomials
$\mathcal{R}$, $ \mathcal{R}_0$ and $\mathcal{R}_{ij}$ $(1 \le i
\le r,\, 1 \le j \le n_i)$ to obtain new polynomials giving the
geometric resolution of the set of common zeros of the
polynomials $F^{(i)}_k$ in $\P^{n_1} \times \dots \times
\P^{n_r}$.

If $(C^{(i)}_{j_1\dots  j_r})_{1\le i \le r,\, 0\le j_t \le n_t}$
are the parameters representing the generic payoffs of a game with
the given structure, in order to obtain the geometric resolution
of the set of its quasi-equilibria it suffices to substitute
$C^{(i)}_{j_1\dots j_{i-1} k j_{i+1} \dots j_r} -
C^{(i)}_{j_1\dots j_{i-1} 0 j_{i+1} \dots j_r}$ for the
coefficient $A^{(i k)}_{J_{-i}}$ for every $1\le i \le r$, $1\le
k\le n_i$ and each $J_{-i} = j_1\dots j_{i-1} j_{i+1} \dots
j_r$ with $0\le j_t \le n_t$ for every $t\ne i$. In this way, the
algorithm obtains polynomials $P(T)$, $ \partial P/\partial T(T)$
and $W_{ij}(T)$ in $\Q[C^{(i)}_{j_1\dots j_r}][T]$ such that the
set of quasi-equilibria of the generic game in $\P^{n_1} \times
\dots \times \P^{n_r}$ is represented as follows:
\begin{equation}\label{quasieq}
\big\{\big((\frac{\partial P}{\partial
T}(t):W_{11}(t):\dots:W_{1n_1}(t)), \dots, (\frac{\partial
P}{\partial T}(t):W_{r1}(t):\dots:W_{rn_r}(t)) \big)\, /\ t \in
\overline{\mathbb{K}}, \ P(t) = 0\big\},
\end{equation}
where $\mathbb{K}:=\Q(C^{(i)}_{j_1\dots j_r})$.  The polynomials
$P$,  $\partial P/\partial T$ and $W_{ij}$ are encoded by
straight-line programs of length $O( D^2 (D + n_1\dots n_r\delta
\log(D) r^2 n^4 (n^3 + r N)))$ over $\Q$, which is also the order
of complexity of the whole computation. The upper bounds  $\deg_T
P\le \delta$ and $\deg_T W_{ij}<\delta$ follow from the stated upper
bounds for the degrees of $\mathcal{R}$ and $\mathcal{R}_{ij}$ in
the variable $A^{(0)}_0$.
\end{proof}

\subsection{Games with the maximum number of totally mixed Nash
equilibria}\label{calculo}

The existence of games with the maximum
possible number of totally mixed Nash equilibria, namely the
multihomogeneous B\' ezout number $\delta$ of the associated
polynomial equation system, was proved in \cite{McKelvey}.
However, no characterization of those games has been provided. In
this subsection, we will obtain a finite family of polynomial
conditions (inequalities over the reals) ensuring
that a given game satisfying those conditions has $\delta$ totally
mixed Nash equilibria. We keep our previous assumptions and
notation.

Let $P(T)$ and $W_{ij}(T)$ $(1\le i \le r, \, 1\le j \le n_i)$ in
$\Q[C^{(i)}_{j_1\dots  j_r}][T]$ be polynomials, as in Theorem
\ref{eq-quasieq}, which give a geometric resolution of the set of
quasi-equilibria (see (\ref{quasieq})) of a game with $r$ players
with $n_1+1, \dots, n_r+1$ strategies and generic payoffs
$(C^{(i)}_{j_1\dots j_r})_{1\le i \le r,\, 0\le j_t \le n_t}$.

Consider a specific choice of payoff values $c:=
(c^{(i)}_{j_1\dots j_r})_{1\le i \le r,\, 0\le j_t \le n_t}$ over
$\R$ and assume that the polynomials $P(c)(T)$ and $W_{ij}(c)(T)$
obtained from  $P(T)$ and $W_{ij}(T)$ by specializing the
parameters at $c$ provide a geometric resolution of the set of
quasi-equilibria of the game with the given payoffs. Then, the
arguments in the last paragraph of subsection \ref{TMNEeq} imply
that the totally mixed Nash equilibria of the game are those
points  $(\xi_1, \dots, \xi_r) \in \R^{n_1 + 1}\times \dots \times
\R^{n_r + 1}$ of the form
$$\xi_i = \Big(\frac{P' (c)(t)}{S_i(c)(t)},
\frac{W_{i1}(c)(t)}{S_i(c)(t)},\dots,\frac{W_{i\,n_i}(c)(t)}{S_i(c)(t)}\Big)
\qquad (1\le i \le r),$$ where $P':= \partial P/ \partial T$ and
$S_i := P' + \sum_{1\le j\le n_i}W_{ij}$, having all their
coordinates real and positive; that is, with $t$ belonging to
$$\{t\in \R : P(c)(t) = 0, \, P'(c)(t)
> 0,\,  W_{ij}(c)(t)
> 0 \ \forall\, 1 \le i \le r, 1 \le j \le n_i\} $$
or $$\{t\in \R : P(c)(t) = 0,\, P'(c)(t) < 0,\, W_{ij}(c)(t) < 0 \
\forall\, 1 \le i \le r, 1 \le j \le n_i\}.$$ Equivalently, $t$
must belong to the intersection
\begin{equation}\label{realpositive}
\bigcap_{1\le i \le r,\, 1\le j \le n_i}\{ t \in \R : P(c)(t) = 0,
\, (P'(c)W_{ij}(c))(t) >0  \}.
\end{equation}

Thus, our polynomial conditions on the payoff vector $c$ which
imply that the associated game has $\delta$ totally mixed Nash
equilibria will state, on the one hand, that the generic
geometric resolution can be specialized into $c$ leading to a
geometric resolution of the set of quasi-equilibria of the game
with the given payoffs and, on the other hand, that the
cardinality of the set introduced in (\ref{realpositive}) equals
$\delta$.

Note that each of the sets appearing in the intersection defined
in (\ref{realpositive}) is described by means of one equality and
one inequality of univariate polynomials over $\R$. In order to
estimate the cardinality  of one of these sets, we will apply the
following well-known result due to Hermite (see, for example,
\cite[Theorem 4.13]{BPRoy}):

\begin{proposition}
Let $p, q\in \R[T]$ be polynomials and assume that $p$ is
square-free. Let $\S_q : \R[T]/(p) \times \R[T]/(p)\to \R$ be the
symmetric bilinear map defined by $\S_q (f,g) = {\rm
Trace}(M_{fqg})$, where for every $H\in \R[T]$, $M_H : \R[T]/(p)
\to \R[T]/(p)$ denotes the linear map $M_H(f) = H \cdot f$. Then
$${\rm Signature}(\S_q) = \# \{\eta \in \R
\,/\, p(\eta) = 0,\, q(\eta) >0\} - \# \{\eta \in \R \,/\,
p(\eta) = 0, \, q(\eta) <0\}.$$
\end{proposition}

The main result of this section is the following:

\begin{theorem}\label{theogeneric} Under the previous assumptions and notations,
there is a family of $n\delta+1$ polynomials $S_0$,
$S_{ij}^{(h)}$, $1\le i \le r,\, 1\le j \le n_i$ and $1\le h \le
\delta$, in $\Q[C^{(i)}_{j_1\dots j_r}]_{1\le i \le r,\, 0\le j_t
\le n_t}\setminus \{ 0 \}$ with total degrees bounded by
$4\delta^2D$, such that for every vector $c:= (c^{(i)}_{j_1\dots
j_r})_{1\le i \le r,\, 0\le j_t \le n_t}$ with real coordinates
satisfying the conditions
$$
S_0(c) \ne 0,\ S_{ij}^{(h)}(c) >0  \quad (1\le i \le r,\, 1\le j
\le n_i, 1\le h \le \delta),
$$
the game with $r$ players with $n_1+1, \dots, n_r+1$ pure
strategies and payoff values given by $c$ has $\delta$ totally
mixed Nash equilibria.

The polynomials $S_0$ and $S_{ij}^{(h)}$ can be obtained
algorithmically within complexity $O(\delta^2(n\delta^2+ L))$
from a straight-line program of length $L$ encoding polynomials
$P, \, W_{ij}$ as in Theorem \ref{eq-quasieq}. The algorithm
computes straight-line programs of length ${\cal
O}(\delta^2(\delta^2 + L))$ which encode these polynomials.
\end{theorem}

\begin{proof} For every $1\le i \le r, \,
1\le k \le n_i$, let $G_k^{(i)} \in \Q[C^{(i)}_{j_1\dots
j_r}][X_1,\dots, X_r]$ be the polynomial obtained after
$F_k^{(i)}$ (see (\ref{genpoly})) by means of the substitution
stated in (\ref{subsgame}). The system $G_k^{(i)}  = 0$ $(1\le i
\le r, \, 1\le k \le n_i)$ defines the quasi-equilibria of a
generic game with $r$ players and $n_1+1,\dots, n_r+1$ pure
strategies respectively.

Let $P(T)$ and $W_{ij}(T)$ $(1\le i \le r, \, 1\le j \le n_i)$ in
$\Q[C^{(i)}_{j_1\dots  j_r}][T]$ be  polynomials as in Theorem
\ref{eq-quasieq} which give a geometric resolution of the set of
common zeros of $G^{(i)}_k$ over an algebraic closure of
$\mathbb{K}:=\Q(C^{(i)}_{j_1\dots  j_r})$.

Denote $P'$ the derivative of the polynomial $P$ with respect to
its main variable $T$. Let us consider the resultant
$$S_0:= \Res_{\delta, \delta - 1}(P(T), P'(T))\in\Q[C^{(i)}_{j_1\dots
j_r}]$$ of $P(T)$ and $P'(T)$ regarded as polynomials in the
single variable $T$ with coefficients in $\Q[C^{(i)}_{j_1\dots
j_r}]$, where the subindices indicate the degrees in the variable
$T$ of $P$ and $P'$ respectively.

Let us observe that, for every real vector $c =(c^{(i)}_{j_1\dots
j_r})_{1\le i \le r,\, 0\le j_t \le n_t}$ with $S_0(c) \ne 0$, the
polynomial $P(c)(T)$ obtained by specializing all the coefficients
of $P(T)$ at $c$ is a non-zero square-free polynomial of degree
$\delta$. Furthermore, the solution set of the system
$G^{(i)}_{k}(c) = 0$, $1\le i \le r, 1\le k \le n_i$, is a
zero-dimensional sub-variety of $\P^{n_1} \times \dots\times
\P^{n_r}$ with $\delta$ distinct points and the polynomials
$P(c)(T), \ W_{ij}(c)(T) \quad (1\le i \le r,\, 1\le j \le n_i)$
give a geometric resolution of this variety.

Therefore, if the condition $S_0(c) \ne 0$ holds, taking into
account that $\deg_T(P(c)(T)) = \delta$, the arguments at the
beginning of this subsection imply that the game with payoff
vector $c$ will have $\delta$ distinct totally mixed Nash
equilibria if and only if every root $t$ of $P(c)$ is real and
satisfies $(P'(c)W_{ij}(c))(t) > 0$ for every $1\le i \le r,\, 1
\le j \le n_i$ (see (\ref{realpositive})). This additional
condition can be restated as
$$ \# \{t\in \R : P(c)(t) = 0,\, (P'W_{ij}(c))(t) > 0\} = \delta \qquad
\forall\, 1\le i \le r,\, 1\le j \le n_i.
$$
Since $\deg_T P = \delta$,  this last condition is equivalent to
\begin{equation}\label{signatura}
{\rm Signature}(\S_{P'W_{ij}(c)}) = \delta \quad \forall\, 1\le i
\le r, 1\le j \le n_i. \end{equation}

Note that for any fixed pair of indices $i, j$, being
$\S_{P'W_{ij}(c)}$  a bilinear form defined over a vector space of
dimension $\delta$, the equality ${\rm Signature}(\S_{P'W_{ij}(c)})
= \delta$ holds if and only if all the eigenvalues of the matrix of
$\S_{P'W_{ij}(c)}$ in an arbitrary basis of the space are positive,
which is in turn equivalent to the fact that the coefficient
sequence of the characteristic polynomial of this matrix, whose
roots are all real, has $\delta$ sign changes. The conditions
stating these sign changes will be the inequalities ensuring that
the game has the maximum number of totally mixed Nash equilibria.

Now, we will obtain the polynomials $S_{ij}^{(h)}$ giving the
inequalities in the statement of the theorem by considering the
bilinear forms $\S_{P'W_{ij}}$ for generic payoffs and performing
all our computations over the parameter field $\mathbb{K} =
\Q(C^{(i)}_{j_1\dots j_r})$, that is, we will consider the
bilinear form $\S_{P'W_{ij}}:\mathbb{K}[T]/(P) \times
\mathbb{K}[T]/(P) \to \mathbb{K}$. Let $\cM^{(ij)}\in
\mathbb{K}^{\delta\times\delta}$ be the matrix of $\S_{P'W_{ij}}$
in the basis $\cB:=\{1, T, \dots, T^{\delta-1}\}$. Note that, for
every real vector $c =(c^{(i)}_{j_1\dots j_r})_{1\le i \le r,\,
0\le j_t \le n_t}$ with $S_0(c) \ne 0$, the coefficient of
$T^\delta$ of the polynomial $P(c)(T)$ is not zero. Then, the
matrix of the bilinear form $\S_{P'W_{ij}(c)}$ can be computed by
specializing $\cM^{(ij)}$ at $c$.

Actually, to deal with polynomials instead of rational functions
of the parameters, we will not compute the matrices $\cM^{(ij)}$
of the bilinear forms $\S_{P'W_{ij}}$ but (scalar) multiples of
them with polynomial entries. In what follows, we will first show
how to obtain these matrices. Then, we will compute the
coefficients of their characteristic polynomials, and we will
estimate their degrees and the overall complexity for their
computation.

Before proceeding, we are going to make some remarks about the
denominator-free computation of traces of multiplication maps in
$\mathbb{K}[T]/(P)$. It is easy to see that for all $H \in
\mathbb{K}[T]$, we have $M_H = H(M_T)$; that is, every
multiplication map is a polynomial function of the multiplication
map $M_T: \mathbb{K}[T]/(P) \to \mathbb{K}[T]/(P)$, whose matrix
in the basis $\cB:=\{1, T, \dots, T^{\delta-1}\}$ can be read off
from $P$. In fact, if we write $P = \sum_{k = 0}^\delta p_k\,
T^k,$ with $p_k \in \Q[C^{(i)}_{j_1\dots j_r}]_{1\le i \le r,\,
0\le j_t \le n_t}$, the matrix of $M_T$ in the basis $\cB$ is
$$ \left( \begin{array}{ccccc} 0 & \dots&\dots &0 &
-\frac{p_0}{p_\delta}\cr 1 & 0 & \dots& \dots& -\frac{p_1}{p_\delta}
\cr 0 & \ddots & \ddots & \cdots& \vdots \cr \vdots& \ddots & \ddots
& 0 & \vdots \cr 0 & \cdots & 0& 1 & -\frac{p_{\delta -1}}{p_\delta}
\end{array} \right).$$
In order not to deal with denominators in our computations in the
case when $H\in \Q[C^{(i)}_{j_1\dots j_r}][T]$, we will consider the
multiplication map $M_{p_\delta T}$, whose matrix in the basis $\cB$
is
$$\cM:=\left(
\begin{array}{ccccc}
0 & \dots &\dots &0 & -p_0\cr p_\delta & 0 &\dots & \dots& -p_1 \cr
0 & \ddots & \ddots & \cdots & \vdots \cr \vdots& \ddots& \ddots & 0
& \vdots \cr 0 & \cdots & 0 & p_\delta & -p_{\delta -1}
\end{array} \right).$$
If $H = \sum_{h = 0}^d b_h T^h$ with $b_h \in
\Q[C^{(i)}_{j_1\dots j_r}]$, we have that ${\rm Trace}(M_H) =
\sum_{h = 0}^d b_h{\rm Trace}((M_T)^h)$ and, therefore, we can obtain 
a multiple of ${\rm Trace}(M_H)$ working over the
polynomial ring $\Q[C^{(i)}_{j_1\dots j_r}]$ due to the following
identity which holds for every $K\ge d$:
\begin{equation}\label{trazas}
p_\delta^K{\rm Trace}(M_H) = \sum_{h = 0}^d p_\delta^{K-h}b_h{\rm
Trace}((M_{p_\delta T})^{h}) .
\end{equation}

Taking into account the previous remark, we are able to obtain a
multiple of the matrix $\cM^{(ij)}$ by a sufficiently big power of
the leading coefficient $p_\delta$ with all its entries in
$\Q[C^{(i)}_{j_1\dots j_r}]$ by working over this polynomial ring:
for every $1\le \alpha, \, \beta \le \delta$, we have that
\begin{equation}\label{matrices}
(\cM^{(ij)})_{\alpha, \beta} = {\rm
Trace}(M_{T^{\alpha-1}P'W_{ij}T^{\beta-1}}) \end{equation} and
then, as $\deg_T(P') <\delta$ and $\deg_T(W_{ij}) < \delta$, the
upper bound $\deg (T^{\alpha-1}P'W_{ij}T^{\beta-1}) \le 4 \delta -
4$ holds. We conclude that the matrix
$p_\delta^{4\delta-4}\cM^{(ij)}$ has all its entries in the
polynomial ring $\Q[C^{(i)}_{j_1\dots j_r}]$. Moreover, as we are
multiplying the matrix $\cM^{(ij)}$ by an \emph{even} power of
$p_\delta$, for every payoff vector $c$ with $S_0(c) \ne 0$, the
identity  $${\rm Signature}(\cM^{(ij)}(c)) = {\rm
Signature}(p_\delta(c)^{4\delta - 4} \cM^{(ij)}(c))$$ holds (note
that $S_0(c) \ne 0$ implies $p_\delta(c) \ne 0$).

Therefore, our algorithm will compute the matrices
$p_\delta^{4\delta-4} \cM^{(ij)}$ and then, their characteristic
polynomials $\cX_{ij} = T^\delta + \sum_{h=1}^{\delta} (-1)^{h}
S_{ij}^{(h)} T^{\delta - h} \in \Q[C^{(i)}_{j_1\dots j_r}][T]$ for
every $1\le i \le r,\, 1\le j \le n_i$. Our previous arguments
imply that for a payoff vector $c$ with $S_0(c) \ne 0$, condition
(\ref{signatura}) is equivalent to
$$S_{ij}^{(h)} (c) > 0 \qquad \forall \, 1\le i \le r,\, 1\le j \le n_i,\, 1\le h \le
\delta.$$

Now we detail the successive steps of the algorithm, we estimate
its complexity and the length of a straight-line program encoding
the output polynomials $S_{ij}^{(h)}$, and we give an upper bound
for the degrees of these polynomials.

The entries of the matrices $p_\delta^{4\delta-4}\cM^{(ij)}$ are
obtained according to identities (\ref{matrices}) and
(\ref{trazas}).

First, from a straight-line program of length $L$ encoding the
polynomials $P$ and $W_{ij}$ $(1\le i \le r,\, 1\le j \le n_i)$,
we obtain a straight-line program of length $O(\delta^2L)$ for the
coefficients of $P$ in the variable $T$ and a straight-line
program of length $O(\delta^2(L+n))$ for the coefficients of all
the polynomials $P'W_{ij}$ $(1\le i \le r, 1\le j \le n_i)$ in the
variable $T$ within complexity $O(\delta^2(L+n))$ (see \cite[Lema
21.25]{Burgisser}).

Then, using the coefficients of $P$, the matrix $\cM \in
\Q[C^{(i)}_{j_1\dots j_r}]^{\delta \times \delta}$ of the
multiplication map $M_{p_\delta T}$ is constructed and the powers
$\cM^h$ for $h = 0, \dots, 4\delta - 4$ are computed. Note that
every matrix in $\Q[C^{(i)}_{j_1\dots j_r}]^{\delta \times
\delta}$ can be multiplied by $\cM$ with ${\cal O}(\delta^2)$
additions and multiplications in $\Q[C^{(i)}_{j_1\dots j_r}]$ and
so, the powers $\cM^h$ for $h\le 4\delta -4$ of this matrix can be
computed within $O(\delta^3)$ polynomial additions and
multiplications. Thus, we obtain straight-line programs of length
$O(\delta^2(\delta + L))$ for the entries of the matrices $\cM^h$
$(0\le h \le 4\delta - 4)$ within complexity $O(\delta^3)$.

Therefore, a straight-line program of length $O(\delta^2(\delta +
L))$ computing ${\rm Trace}((M_{p_\delta T})^h)$, for every $0\le
h\le 4\delta - 4$ is obtained with $O(\delta^2)$ additional
operations.

Now, fix $i, j$ with $1\le i \le r$ and $1\le j \le n_i$. Let
$P'W_{ij} = \sum_{h=0}^{2\delta - 2} b^{(ij)}_h T^h$ be the
expansion of $P'W_{ij}$ in the variable $T$ with coefficients in
$\Q[C^{(i)}_{j_1\dots j_r}]$. Recall that in the first step of the
algorithm we have obtained a straight-line program of length
$O(\delta^2 (L+n))$ for the polynomials $b^{(ij)}_h$. For every
$\alpha, \beta$ with $1\le \alpha, \beta \le \delta$, we obtain a
straight-line program for $(p_\delta^{4\delta-4}
\cM^{(ij)})_{\alpha, \beta}$ by means of the formula:
$$p_{\delta}^{4\delta-4}{\rm Trace}(M_{T^{\alpha-1}P'W_{ij}T^{\beta-1}})
 =
\sum_{h = \alpha + \beta-2}^{2\delta + \alpha + \beta-4}
p_\delta^{4\delta-4-h}b^{(ij)}_{h-\alpha- \beta+2}{\rm
Trace}(M_{(p_\delta T)^h})$$ within additional complexity
$O(\delta)$. The length of this straight-line program is $
O(\delta^2(\delta + L+n))$.

The complexity of the computation of the entries of all the
matrices $p_\delta^{4\delta-4} \cM^{(ij)}$ from the coefficients
$b^{(ij)}_h$ and the traces  ${\rm Trace}((M_{p_\delta T})^h)$ is
$O(n\delta^2)$.

Finally, we apply the division-free algorithm described in
\cite{Berk} in order to compute the coefficients of the
characteristic polynomial $\cX_{ij}$ of each of these matrices.
For every $1\le i \le r, 1 \le j \le n_i$, the algorithm produces
straight-line programs of length $O(\delta^2(\delta^2 +L+n))$
encoding these coefficients within complexity $O(\delta^4)$.

The overall complexity of the algorithm is $O(\delta^2(n\delta^2 +
L))$.

Let us observe that the entries of the matrices
$p_\delta^{4\delta-4} \cM^{(ij)}$ are polynomials in
$\Q[C^{(i)}_{j_1\dots j_r}]$ of total degrees bounded by $4\delta
D$ and so, the polynomials $S_{ij}^{(h)}$ obtained by the
described procedure, which are the (signed) coefficients of the
characteristic polynomials $\cX_{ij}$, have total degrees bounded
by $4 \delta^2  D$.
\end{proof}

Note that, with the same notation as in the previous theorem,
for a generic game with the given structure (namely, any game
whose payoff vector $c$ satisfies $S_0(c)\ne 0$) the conditions
$S^{(h)}_{ij}(c)>0$ are \emph{equivalent} to the fact that the considered game
has the maximum possible number of totally mixed Nash equilibria.

\section{The set of totally mixed Nash equilibria of an arbitrary game}\label{sec particular}

When dealing with a particular game with specific payoff values,
the geometric resolution of the set of quasi-equilibria of a
generic game with the same structure might not be useful in order
to describe the quasi-equilibria of the given game. For instance,
the minimal polynomial $P(T)$ or its derivative $\frac{\partial
P}{\partial T}(T)$ computed by the algorithm shown in Theorem
\ref{eq-quasieq} could vanish identically when specialized at the
considered payoff values. This section is aimed at adapting the
procedures previously developed in order to handle these
particular situations.

First, we show an algorithm computing a geometric resolution of
the set of affine quasi-equilibria of the game provided that it
has a zero-dimensional set of quasi-equilibria. In a second step,
this algorithm is extended to find all \emph{isolated} affine
quasi-equilibria of an arbitrary game.

Once a finite set of affine quasi-equilibria of a game including
the isolated ones is computed, we apply standard algorithms of
semialgebraic geometry to estimate the number of isolated totally
mixed Nash equilibria of the game.

\subsection{Games with a zero-dimensional set of quasi-equilibria}

As in the previous sections, we consider an $r$-person
non-cooperative game in normal form in which the players have
$n_1+1,\dots,n_r+1$ distinct available pure strategies each. Let
$c^{(i)}:= (c^{(i)}_{j_1,\dots,j_r})_{0\le j_k\le n_k}$ denote the
payoff matrix to player $i$ for every $1\le i \le r.$

Then, the polynomial equations defining the Nash equilibria of the
game (see Section \ref{TMNEeq}) can be obtained by specializing
the coefficients of each of the generic multilinear polynomials
$F^{(i)}_k$ introduced in (\ref{genpoly}) in the vector $a^{(ik)}
:= ( a^{(ik)}_{j_1 \dots j_{i-1} j_{i+1} \dots j_r})_{0\le
j_t \le n_t}$ defined from the payoffs of the game as follows
$a^{(ik)}_{j_1 \dots j_{i-1} j_{i+1} \dots j_r}: =
c^{(i)}_{j_1\dots j_{i-1} k j_{i+1} \dots j_r} - c^{(i)}_{j_1\dots
j_{i-1} 0 j_{i+1} \dots j_r}$. Thus, if $a:= (a^{(ik)})_{1\le
i\le r,\, 1\le k\le n_i}$, the set of quasi-equilibria of the game
is the algebraic variety
$$V_a := \{ \xi:= (\xi_1,
\dots, \xi_r)\in \P^{n_1} \times \dots \times \P^{n_r} \ / \
F^{(i)}_k(a, \xi) = 0 \ \forall 1\le i \le r,\, 1\le k \le n_i
\}.$$

In this section we consider the case when $V_a$ is
zero-dimensional. In order to decide whether this is the case for
given payoff values, we study the non vanishing of an adequate
multihomogeneous resultant. More precisely, we consider a generic
polynomial of multi-degree $d_0:=(1,\dots,1)$ in the groups of
variables $X_1,\dots,X_r$,
$$F_{0} = \displaystyle\sum_{\tiny\begin{array}{c}  1\le i \le r\\ 0\le j_i\le n_i\end{array}}
 A^{(0)}_{j_1\dots j_r}
x_{1j_1} \dots x_{rj_r}$$ and the resultant  $R := {\rm Res}(F_0,
F_1^{(1)}, \dots, F_{n_1}^{(1)}, \dots, F_1^{(r)}, \dots,
F_{n_r}^{(r)})$. Let $R_a(A^{(0)}_{j_1\dots j_r})$ be the
polynomial obtained by substituting the coordinates of $a=
(a^{(ik)})_{1\le i\le r,\, 1\le k\le n_i}$ for the variables
$A^{(ik)}$ in the polynomial $R$. Then:

\begin{remark}
The algebraic variety $V_a$ is either zero-dimensional or empty if
and only if the polynomial $R_a$ is not identically zero.
\end{remark}

\begin{proof}{Proof.} If $V_a$ is empty, the result is straightforward. If $V_a$ is
zero-dimensional, there exists a multilinear polynomial $f_0\in
\Q[X_1,\dots, X_r]$ which does not vanish at any of the (finitely
many) points of $V_a$ and therefore, $R_a$ does not vanish at the
coefficients of $f_0$.

On the other hand, if  the projective variety $V_a$ has positive
dimension, any multilinear polynomial $f_0$ has zeros in $V_a$.
Therefore, $R_a$ is identically zero.
\end{proof}

Now we will show how to decide whether $R_a$ is identically zero
algorithmically and we will estimate the complexity of this
procedure.

First, we compute the multihomogeneous resultant $R$ by means of
the algorithm described in \cite{JeSa} adapted according to
Subsection {\ref{appendix}} below. Thus, a straight-line program
of length $\mathcal{L}:=O( D^2 (D + n_1\dots n_r\delta \log(D) r^2
n^4 (n^3 + r N)))$ encoding $R$ is obtained within complexity of
the same order as $\mathcal{L}$. The specialization
$A^{(ik)}:=a^{(ik)}$ which leads us to the polynomial $R_a$ does
not modify this complexity order.

Since $R_a$ is a multivariate polynomial in $\prod_{1\le i \le r}
(n_i+1)$ variables encoded by a straight-line program, we do not
know its coefficients; moreover, the complexity to compute them is
exponential. Therefore, in order to decide whether it is the zero
polynomial or not, we will make a suitable specialization so that
solving this problem amounts to solving the same problem in the
univariate setting.

More precisely, we specialize the variable coefficients
$A^{(0)}_{j_1\dots j_r}$ into successive powers of a new variable
$t$, that is,
\begin{equation} \label{especializacion} A^{(0)}_{j_1\dots j_r} =
t^{j_1 + (n_1 + 1)j_2 + (n_1 + 1)(n_2 + 1)j_3 + \dots +
(\prod_{j=0}^{r-1}(n_j + 1))j_r }.\end{equation} Making this
specialization in $F_0$, we obtain a new polynomial $f_0 \in
\mathbb{Q}[t][X_1,\dots, X_r]$. Let us show that if $V_a$ consists
of finitely many points, there is a polynomial $f_0(t_0)$ with
coefficients of this type (namely, successive powers of the same
scalar $t_0$) which does not vanish at any point of $V_a$.

For every $\xi:= (\xi_1, \dots, \xi_r)\in \P^{n_1} \times \dots
\times \P^{n_r}$, we have $$f_0(\xi)(t) = \sum_{j_1\dots j_r}
\xi_{1j_1}\dots\xi_{rj_r}t^{j_1 + (n_1 + 1)j_2 + (n_1 + 1)(n_2 +
1)j_3 + \dots + (\prod_{j=0}^{r-1}(n_j + 1))j_r} \in \C[t],$$
which is a non-zero polynomial, due to the fact that there exists
at least one choice of $j_1,\dots,j_r$ for which the product
$\xi_{1j_1}\dots\xi_{rj_r}$ is not $0$. Now, if $V_a$ is a finite
set, we may consider the polynomial $\Delta(t):= \prod_{\xi \in
V_a} f_0(\xi)(t)$. Note that $\Delta \in \C[t]$ is a non-zero
polynomial and, therefore, there exists $t_0 \in \mathbb{Q}$ with
$\Delta(t_0)\ne 0$. Due to the definition of $\Delta$, the
polynomial $f_0(t_0)\in \mathbb{Q}[X_1,\dots,X_r]$ does not vanish
at any point of $V_a$. This implies that the resultant of
$f_0(t_0), F^{(i)}_k(a^{(ik)})$ $(1\le i \le r,\, 1\le k \le n_i)$
is not zero, that is, the polynomial $R_{a,t}\in \mathbb{Q}[t]$
obtained by specializing $R_a$ following (\ref{especializacion})
does not vanish at $t_0$. We conclude that $R_{a,t}$ is not the
zero polynomial.

Therefore, in order to decide whether $R_a\in
\mathbb{Q}[A^{(0)}_{j_1\dots j_r}]$ is zero or not, it suffices to
decide whether the associated polynomial $R_{a, t}\in \mathbb{Q}[t]$
is zero or not. Taking into account that $\deg(R_{a, t}) \le
(\prod_{1\le i \le r}(n_i+1)-1)\delta$, it follows that this task
can be achieved by evaluating $R_{a, t}$ at $(\prod_{1\le i \le
r}(n_i+1))\delta$ distinct values of $t$, which is done by
substituting the powers of these values for  the variables
$A^{(0)}_{j_1\dots j_r}$ according to (\ref{especializacion}) in the
straight-line program for $R_a$.

 The overall complexity of this procedure is
$\Big((\prod_{1\le i \le r}(n_i+1))+ {\cal L} \Big)(\prod_{1\le i
\le r}(n_i+1))\delta = O( D^2 (n_1\dots n_r)^2 \delta (D +
n_1\dots n_r\delta \log(D) r^2 n^4 (n^3 + r N))) $.

\bigskip

Assume now that the set of quasi-equilibria of the given game is
finite, that is, the polynomial $R_a$ is not identically zero.

\begin{proposition}\label{affquasi0}
Following the previous assumptions and notation, there is an
algorithm which computes a geometric resolution of the set of
affine quasi-equilibria of a game with $r$ players having
$n_1+1,\dots, n_r+1$ pure strategies respectively within
complexity $O(\delta^8D^2 (D +n_1\dots n_r\delta \log(D) r^2 n^5
(n^3 + r N)))$ provided that the associated set of
quasi-equilibria $V_a$ is zero-dimensional.
\end{proposition}

\begin{proof}{Proof.} We will obtain the desired geometric resolution
from the specialized resultant $R_a$. In order to do this, we first
show that this specialized resultant factorizes as follows:
$$
R_a = C_a\prod_{\xi \in V_a} F_0(\xi)^{m_\xi}, \ \hbox{ with } C_a
\in \Q,\ m_\xi\in \N.
$$
Set $\cF_a := \prod_{\xi\in V_a}F_0(\xi)\in \Q[A^{(0)}]$. For a
given coefficient vector $a^{(0)}$, we have that $\cF_a(a^{(0)}) =
0$ if and only if there is a point  $\xi \in V_a$ such that
$F_0(a^{(0)},\xi) = 0$; that is, if and only if $F_0(a^{(0)})$,
$F^{(1)}_{1}(a^{(11)}), \dots, F^{(r)}_{n_r}(a^{(rn_r)})$ have a
common root in  $\P^{n_1} \times \dots \times \P^{n_r}$. But this
last condition is equivalent to the fact that the resultant ${\rm
Res}(F_0, F^{(1)}_1,\dots, F^{(r)}_{n_r})$ vanishes at
$a^{(0)},a^{(11)}, \dots, a^{(rn_r)}$ or, equivalently, that
$R_a(a^{(0)}) = 0$. We conclude that the polynomials  $\cF_a$ and
$R_a$ have the same zero set and so, they have the same
irreducible factors, which proves our assertion.

We are interested in describing the set of \emph{affine}
quasi-equilibria of the game, that is, the variety $V_a^{{\rm
aff}} := \{ \xi \in V_a : \xi_{i0} \ne 0 \ \forall\, 1\le i \le
r\}$. For this reason, we are going to compute the polynomial
$R_a^{{\rm aff}}:= \prod_{\xi \in V_a^{{\rm aff}}}
F_0(\xi)^{m_\xi}$. Note that
\begin{equation}\label{forma}
R_a = C_a \prod_{\xi \in V_a-V_a^{{\rm aff}}} F_0(\xi)^{m_\xi}
\prod_{\xi \in V_a^{{\rm aff}}} F_0(\xi)^{m_\xi}, \end{equation}
where $R_a^{{\rm aff}}=\prod_{\xi \in V_a^{{\rm aff}}}
F_0(\xi)^{m_\xi}$ is monic in the variable $A^{(0)}_{0\dots 0}$
and $\prod_{\xi \in V_a-V_a^{{\rm aff}}} F_0(\xi)^{m_\xi}$ does
not depend on this variable. Then, $C_a \prod_{\xi \in
V_a-V_a^{{\rm aff}}} F_0(\xi)^{m_\xi}$ is the leading coefficient
of  $R_a$ in the variable $A^{(0)}_{0\dots 0}$. Therefore,
$R_a^{{\rm aff}}$ can be obtained from $R_a$ by dividing it by
this leading coefficient.

Let us consider the generic form $f_0:= A^{(0)}_{0} x_{10} \dots
x_{r0} +
\displaystyle\sum_{\tiny\begin{array}{c} 1\le i\le r\\
1\le j \le n_i\end{array}} \hspace{-5mm} A^{(0)}_{ij} x_{10} \dots
x_{i-1 \, 0}x_{i j}x_{i+1 \, 0} \dots x_{r0}$.  After specializing
$R_a^{{\rm aff}}$ as follows:
\begin{equation}\label{especializaciondos}\begin{array}{rcll}
A^{(0)}_{0\dots 0} & \mapsto & A^{(0)}_0\\
A^{(0)}_{0\dots0\, j\, 0\dots 0} & \mapsto & A^{(0)}_{ij} &
\hbox{ where the index $j$ is in the $i$th place} \ (1\le j \le
n_i)\cr A^{(0)}_{j_1\dots j_r} &\mapsto& 0 & \hbox{ otherwise
}\end{array}
\end{equation} we obtain
the polynomial $\widetilde P_a:= \prod_{\xi \in V_a^{{\rm aff}}}
f_0(\xi)^{m_\xi}$. The geometric resolution of the variety
 $V_a^{\rm aff}$ can be obtained from the square-free part  $P_a:= \prod_{\xi \in V_a^{{\rm aff}}}
f_0(\xi)$ of $\widetilde P_a$ and the partial derivatives of $P_a$
as shown in the proof of Theorem \ref{eq-quasieq} and substituting
the powers of a conveniently chosen scalar for the variables
$A_{ij}$.

Now we estimate the complexity of the different steps of our
algorithm.

Computation of $R_a^{{\rm aff}}$: First, we compute the degree
$d_a:=\deg_{A^{(0)}_{0\dots 0}}(R_a)$. To do this, let $t_0 \in
\Q$ be the element we have already obtained such that $R_{a, t}
(t_0) \ne 0$. Let $R_a^{(t_0)} \in \Q[A^{(0)}_{0\dots 0}]$ be
the (non-zero) polynomial obtained from $R_a$ after specializing
it as in (\ref{especializacion}) for every $(j_1,\dots, j_r) \ne
(0,\dots,0)$ setting $t=t_0$. Because of (\ref{forma}),
$\deg(R_a^{(t_0)}) = d_a$. Then, in order to compute $d_a$, it
suffices to compute the coefficients of $R_a^{(t_0)}$ up to degree
$\delta$ (which is an a priori upper bound for $d_a$). The
complexity of this computation is of order $O(\delta^2{\cal L})$.
Now we obtain an slp of length $O(\delta^2 {\cal L })$ for the
coefficient of $(A^{(0)}_{0\dots 0})^{d_a}$ of $R_a$ and we
obtain $R_a^{{\rm aff}}$ by dividing $R_a$ by this coefficient.
As the divisor does not vanish when its variables are specialized
in the successive powers of $t_0$, this polynomial division can be
done within complexity $O(\delta^4{\cal L})$, and produces an slp
of the same order (\cite{Strassen}).

Computation of the polynomial $P_a$: We obtain $P_a$ as the
quotient $ \widetilde P_a / {\rm gcd}(\widetilde P_a,
\frac{\partial{\widetilde P_a}}{\partial A^{(0)}_0})$. The
specialization in (\ref{especializaciondos}) gives an slp for
$\widetilde P_a$ without modifying the previous length.

To obtain the square-free part of $\widetilde P_a$, we apply a
well-known subresultant-based procedure for the computation of the
gcd of two polynomials (see, for instance, \cite{BT}). We proceed
as before by specializing all the variables in $\widetilde P_a$,
except for the main variable $A^{(0)}_0$ into the successive
powers of a new variable $t$ and deal with the polynomial we
obtain in $\Q[t][A^{(0)}_0]$. We will use the following notation:
given $F\in \Q[A^{(0)}_0,A^{(0)}_{ij}]$,  $F^{(t)}$ will denote
the polynomial obtained from $F$ by doing the previously mentioned
specialization.

Let $G = {\rm gcd}(\widetilde P_a, \frac{\partial{\widetilde
P_a}}{\partial A^{(0)}_0})$. As, for $\xi_1 \ne \xi_2$ we have
that $f_0^{(t)}(\xi_1)$ y $f_0^{(t)}(\xi_2)$ are relatively prime
irreducible polynomials in $\overline \Q[A^{(0)}_0, t]$, we have
that
$$G^{(t)}  = \prod_{\xi\in V_a^{{\rm aff}}}  f_0^{(t)}(\xi)^{{m_\xi}-1} =
{\rm gcd}(\widetilde P_a^{(t)}, \frac{\partial{\widetilde
P_a^{(t)}}}{\partial A^{(0)}_0}).$$ Then, we can obtain $\tilde
d_a := \deg({\rm gcd}(\widetilde P_a, \frac{\partial{\widetilde
P_a}}{\partial A^{(0)}_0}))$ using the polynomials $\widetilde
P_a^{(t)}$ and $\frac{\partial{\widetilde P_a^{(t)}}}{\partial
A^{(0)}_0}$. In order to compute the degree of the gcd of these
polynomials, we look for their first non-zero subresultant. In
each step, to decide whether the considered subresultant (which is
a polynomial of degree at most $2\delta^2n$ in $\Q[t]$) is zero or
not, we evaluate the variable $t$ in a sufficient number of
elements of $\Q$. First, we obtain an slp of length $O(\delta^4
\cL +n)$ for $\widetilde P_a^{(t)}$ and then an slp of length
$O(\delta^2 (\delta^4 \cL +n))$ for its coefficients in the
variable $A^{(0)}_0$. For a specific evaluation of $t$ the
complexity of the computation of all the subresultants is
$O(\delta^6{\cal L}+ \delta^2 n)$, and therefore, the whole
complexity of this step is bounded by $O(\delta^8{\cal L}n +
\delta^4 n^2)$. We obtain $\widetilde G = {\rm sr}_a \, . \,{\rm
gcd}(\widetilde P_a, \frac{\partial{\widetilde P_a}}{\partial
A^{(0)}_0})$  as the polynomial subresultant of index $\tilde
d_a$, where ${\rm sr}_a$ is the main subresultant of index $\tilde
d_a$, within complexity $O(\delta^6{\cal L})$.

Finally, the polynomial $P_a$ is obtained by dividing  ${\rm
sr}_a\widetilde P_a$ by $\widetilde G$. Note that we already know
a point $\tau$ where ${\rm sr}_a$ does not vanish ($\tau$ is the
value obtained in the previous step when computing $\tilde d_a$).
Then, if we consider the nonzero polynomial $ \widetilde
G^{(t)}(A^{(0)}_0, \tau) \in \Q[A^{(0)}_0]$ of degree $\tilde
d_a$, by evaluating it in at most $\tilde d_a+1$ elements in $\Q$,
we obtain  $a^{(0)}_0\in \Q$ such that $\widetilde
G^{(t)}(a^{(0)}_0, \tau) \ne 0$. This enables us to compute the
quotient $P_a$ by applying the classical division avoiding
algorithm in \cite{Strassen}. The complexity of this step is of
order $O(\delta^8{\cal L})$.

Observe that the linear form $\ell$ whose coefficients are $(-
\tau^k)_{0\le k \le n-1}$ is a separating linear form for
$V_a^{{\rm aff}}$. To compute the required geometric resolution
associated with $\ell$ from the polynomial $P_a$ we proceed as in
the proof of Theorem \ref{eq-quasieq}: we specialize $P_a$ and its
partial derivatives at the coefficients of $\ell$.

The overall complexity of the algorithm  is $O(\delta^8{\cal L}n
+\delta^4 n^2) = O(\delta^8D^2 (D + n_1\dots n_r\delta \log(D) r^2
n^5 (n^3 + r N)))$.
\end{proof}

 Now, we are able to compute the number of totally mixed
Nash equilibria of a game with zero-dimensional set of
quasi-equilibria from the geometric resolution given by
Proposition \ref{affquasi0}.

Let $p$ and $w_{ij}$ $(1\le i \le r,\, 1\le j \le n_i)$ be the
polynomials giving the geometric resolution of $V_a^{\rm aff}$.
Then, the totally mixed Nash equilibria of the game are the points
$(\xi_1,\dots, \xi_r)$ with $\xi_i=(p'(t)/s_i(t),
w_{i1}(t)/s_i(t),\dots, w_{ir}(t)/s_i(t))$, where $s_i=
p'+\sum_{j=1}^{n_i} w_{ij}$ for $i=1,\dots, r$, and $t$ is a root of
$p$, having all their coordinates real and positive.

Then, the number of totally mixed Nash equilibria of the game equals
the cardinality of the union of the sets
$$\{ t\in \R : p(t) = 0, p'(t)>0,\, w_{ij}(t)>0 \ \forall\, 1\le i \le
r, \, 1\le j \le n_i\}$$ and
$$\{ t\in \R : p(t) = 0, p'(t)<0,\, w_{ij}(t)<0 \ \forall\, 1\le i \le
r, \, 1\le j \le n_i\}.$$

Once we have the polynomials $p$, $p'$ and $w_{ij}$ $(1\le i \le
r, \, 1\le j \le n_i\})$ encoded in dense form, by using the
algorithm in \cite[Section 3.3]{Canny}, it is possible to compute
this cardinality within  complexity $O(n\delta^3)$.

Therefore we have:

\begin{proposition}
Following the previous assumptions and notation, there is an
algorithm which computes the number of totally mixed Nash
equilibria of a game with $r$ players having $n_1+1,\dots, n_r+1$
pure strategies respectively within complexity $O(\delta^8D^2 (D +
n_1\dots n_r\delta \log(D) r^2 n^5 (n^3 + r N)))$  provided that
the associated set of quasi-equilibria of the game is
zero-dimensional.
\end{proposition}

\subsection{Computing the isolated affine quasi-equilibria of an arbitrary
game}

When we are dealing with an arbitrary game, it may happen that the
set $V_a$ of its quasi-equilibria has positive dimension and,
therefore, the polynomial $R_a$ introduced in the previous
subsection is identically zero. In this case, we are going to use
a procedure applied in \cite{HJSS05} to obtain a non-zero multiple
of the minimal polynomial of a generic linear form over the
isolated points of the set $V_a^{{\rm aff}}$ of affine
quasi-equilibria. For the sake of completeness, we are going to
explain this procedure briefly.

We consider a \emph{sufficiently generic} coefficient vector 
$b:=(b^{(ik)})_{1\le i \le r, 1\le k \le n_i}$
such that $R_b \not\equiv 0$ (we remark that the vector $b$ can be
chosen at random or effectively constructed as the coefficient
vector a system with $\delta$ many common roots). Then, if $\tau$
is a new variable, the polynomial $R_{a+ \tau(b-a)}$ is non-zero.
If $f^{(i)}_k$ denotes the polynomial obtained from $F^{(i)}_k$ by
evaluating $x_{j0} = 1$ $(1\le j \le r)$, let $\widetilde
f^{(i)}_k:= f^{(i)}_k(a^{(ik)} +\tau (b^{(ik)}- a^{(ik)}))$
$(1\le i \le r,\, 1\le k \le n_i)$. Let $\widetilde
P_{a+\tau(b-a)}$ be the polynomial obtained from $R_{a+
\tau(b-a)}$ by specializing it as in (\ref{especializaciondos})
and let $L:= \sum_{i,j} A_{ij}^{(0)}x_{ij}$.

As $\widetilde P_{a+\tau(b-a)}(\tau, A^{(0)})\mid_{\tau=0} \equiv
0$, there exists $\epsilon>0$ such that $\widetilde
P_{a+\tau(b-a)}= \tau^\epsilon \widetilde P$ for a polynomial
$\widetilde P \in \Q[\tau, A^{(0)}]$ with $\widetilde
P(\tau,A^{(0)})\mid _{\tau=0}\not\equiv 0$. Let $P:=  \widetilde
P\mid_{\tau=0, A_0^{(0)} = - T}$. We are going to show that this
polynomial $P$ is a multiple of the minimal polynomial of $L$ over
the set of isolated points of $V_a^{{\rm aff}}$.

Since the polynomial $R$ is a linear combination of $F_0,
F^{(i)}_k$ $(1\le i \le r,\, 1\le k \le n_i)$, we have that $
\widetilde P_{a+\tau(b-a)}\mid_{A_0^{(0)} = - L} = \tau^\epsilon
\widetilde P \mid_{A_0^{(0)} = - L} \in(\widetilde f^{(i)}_k: 1\le
i \le r,\, 1\le k \le n_i)\in \Q[\tau, x_{ij}: 1\le i \le r, 1\le
j\le n_i]$.

Now, taking into account that the number of variables $\tau,
x_{ij}$ $(1\le i \le r, 1\le j \le n_i)$ is greater than the number
of generators of the ideal $(\widetilde f^{(i)}_k: 1\le i \le r,\,
1\le k \le n_i)$, we deduce that each irreducible component of the
variety $V_\tau$ defined by this ideal in $\A^{n+1}$ has dimension
at least $1$. Then, for each isolated point $\xi$ of $V_a^{{\rm
aff}}$, the point $(0, \xi) \in \A^{n+1}$ lies in an irreducible
component $C$ of $V_\tau$ such that $\tau\notin I(C)$. Therefore
$\widetilde P \mid_{A_0^{(0)} = - L} \in I(C)$ and $P\mid_{T=L}$
vanishes at $\xi$.

Now, we are going to show how to compute the polynomial $P$.

Let $\widetilde P_{a+ \tau(b-a)} = \sum_{j=0}^\delta p_j(\tau)
(A^{(0)}_{0})^j$, where $p_j(\tau)$ are polynomials in
$\Q[A^{(0)}_{ij}][\tau]$ such that $p_j(\tau) \not\equiv 0$ for some
$0\le j\le D$. Then $\epsilon := \max\{ k : \tau^k \hbox{ divides }
p_j \ \forall\, 0\le j \le D\}$.

The procedure for the computation of $P$ from a straight-line
program encoding $R$ runs as follows:

\begin{itemize}
\item[(1)] Obtain an slp for $\widetilde P_{a+\tau(b-a)}$.
\item[(2)] Obtain an slp encoding the coefficients $p_j(\tau)$ $(0\le j \le \delta)$.
\item[(3)] For $j=0, \dots \delta$:
\begin{itemize}
\item[(3.a)] Obtain an slp for the coefficients of $p_j(\tau) =
\sum_{k=0}^D p_{jk}\, \tau^{k}$.
\item[(3.b)] Compute $\epsilon_j:= \min\{ k : p_{jk} \ne 0\}$.
\end{itemize}
Compute $\epsilon := \min \{ \epsilon_j : 0\le j \le \delta\}$.
\item[(4)] Obtain an slp for $ P= \sum_{j=0}^\delta
(-1)^jp_{j\epsilon} T^j$ from the slp's encoding $p_{j\epsilon}$
$(0\le j \le \delta)$.
\end{itemize}

Now, we are going to estimate the complexity of this procedure. Let
$\mathcal{L}$ be the length of a straight-line program encoding $R$.

In step (1), we compute first $a+\tau(b-a)$ with $3N$ operations
and then, we specialize the slp for $R$ according to
(\ref{especializaciondos}). We obtain an slp of length
$3N+\mathcal{L}$ for $\widetilde P_{a+\tau(b-a)}$. The
interpolation in step (2) takes $O(\delta^2 (N+\mathcal{L}))$
operations. The complexity of step (3.a) is of order
$O(D^2\delta^2 (N+\mathcal{L}))$.  In order to determine whether
each of the multivariate polynomials $p_{jk} \in \Q[A^{(0)}_{ij}]$
encoded by an slp is zero or not, we apply the probabilistic
Zippel-Schwartz zero test (see \cite{Zippel}). The overall
complexity of step (3) is $O(D^2\delta^3 (N+\mathcal{L}))$. Step
(4) does not modify the order of the complexity.

{}From the polynomial $P$ we have computed we can obtain a geometric
resolution for a finite set of points including the isolated points
of $V_a^{{\rm aff}}$ in the same way we showed in the proof of
Theorem \ref{eq-quasieq}, obtaining a description for the isolated
affine quasi-equilibria of the game. Therefore, we have:

\begin{theorem}
Let notations be as in Theorem \ref{eq-quasieq}. There is a
probabilistic algorithm which computes, given a game with $r$
players having $n_1 +1, \dots, n_r+1$ pure strategies
respectively, a geometric resolution of a finite set of points
including the isolated affine quasi-equilibria of the game within
complexity
 $O( D^4 \delta^3 (D + n_1\dots n_r\,\delta \log(D)
r^2 n^4 (n^3 + r N)))$.
\end{theorem}

Proceeding as in the previous section, from the geometric resolution
given by this theorem, we are able to obtain an upper bound on the
number of isolated  (in the complex space) totally mixed Nash
equilibria of the game within the same order of complexity as in the
previous theorem.

\section{Resultants in the multihomogeneous setting}

\subsection{Computing resultants}\label{appendix}

In this subsection we are going to show how the algorithm in
\cite{JeSa} can be adapted in order to compute the
multihomogeneous resultants we use. That procedure computes
multihomogeneous resultants under the assumption that the
coordinates of each multidegree are all positive and this is not
the case in our situation. 

First, we will prove that the multihomogeneous resultant we want to 
compute is a non-constant polynomial and then, we will show how to modify 
the Poisson formula on which the algorithm in
\cite{JeSa} relies  All the computations required once the Poisson
formula is recursively applied run in the same way as in \cite{JeSa}.

In order to do this, we are going to use the theory in
\cite{Stu94} and \cite{Min}. We can apply these results to our situation because the
multihomogeneous resultant of a family of multihomogeneous
polynomials $G_0,\dots, G_n$ in $r$ groups of variables
$X_j=(x_{j0}, \dots, x_{j n_j})$, with $\sum_{j=1}^r n_j =n$,
coincides with the sparse resultant of the dehomogenized
polynomials $g_0,\dots, g_n$ obtained by setting $x_{j0}=1$ for
every $1\le j \le r$.

Let $\cA_0,\dots, \cA_n\subset \Z^n$ be finite sets and let
$g_0,\dots, g_n$ be polynomials with supports $\cA_0, \dots, \cA_n$
respectively.

For any subset $J\subseteq \{0,\dots, n\}$, let $\cL_J$ be the
lattice generated by $\sum_{j\in J} \cA_j$. Following
\cite{Stu94}, if $I\subset \{0,\dots, n\}$, the collection of
supports $\{\cA_i\}_{i\in I}$ is said to be \emph{essential} if
${\rm rank}(\cL_I) = \# I-1$ and ${\rm rank}(\cL_J) \ge \# J$ for
each proper subset $J$ of $I$.

If there is a unique subcollection $\{\cA_i\}_{i\in I}$ which is
essential, the resultant ${\rm Res}(g_0,\dots, g_n)$ coincides
with the resultant ${\rm Res}(g_i; i \in I)$.

Note that if $G$ is a multihomogeneous polynomial in $r$ groups of
$n_1+1,\dots, n_r+1$ variables each with multidegree
$(\nu_1,\dots, \nu_r)$, the set of exponent vectors of the
dehomogenized polynomial $g$ is $\Delta_{n_1, \nu_1}\times \dots
\times \Delta_{n_r, \nu_r}\subset \Z^{n_1}\times \dots\times
\Z^{n_r},$ where $\Delta_{n_i,\nu_i}$ is $\nu_i$ times the unitary
simplex in $\Z^{n_i}$ for $i=1,\dots, r$.

In our particular case, we will deal with families of
multihomogeneous polynomials in $r$ groups of $m_1+1,\dots, m_r+1$
variables respectively and multidegrees $(1,\dots, 1)$ or
$d_i=(1,\dots, 0,\dots, 1)$ (with $0$ in the $i$th coordinate) for
$i=1,\dots, r$. If one of the polynomials in the family indexed by
$I$ has multidegree $(1,\dots, 1)$, then ${\rm rank}(\cL_I) =
\sum_{i=1}^r m_i$, and the same holds if the family contains
polynomials with multidegrees $d_i$ and $d_j$ with $i\ne j$. If
all the polynomials in the family have multidegree $d_i$, then
${\rm rank}(\cL_I) = \sum_{j\ne i } m_j$.

The resultant will be computed recursively by applying Poisson's
formula (\cite[Lemma 13]{Min}). At each step, we will have to
compute a multihomogeneous resultant in one of the following
settings, where $m_1,\dots, m_r\in \N$ and $m:= m_1+\dots+m_r$:
\begin{itemize}
\item[(I)] $m+1$ multihomogeneous polynomials in $r$ groups of
$m_1+1,\dots, m_r+1$ variables: one polynomial with multidegree
$(1,\dots, 1)$ and, for every $1\le i \le r$, $m_i$ polynomials
with multidegrees $d_i:=(1,\dots, 0,\dots, 1)$ with $0$ in the
coordinate $i$, under the assumption that $m_i \le \sum_{j\ne i}
m_j$ for every $1\le i \le r$, \item[(II)] $m+1$ multihomogeneous
polynomials with multidegrees $(1,\dots, 1)$ in $r$ groups of
$m_1+1,\dots,m_r+1$ variables each, \item[(III)] $m$
multihomogeneous polynomials in $r$ groups of variables with
$m_1,m_2+1,\dots, m_r+1$ variables each, with $m_i$ polynomials of
multidegree $d_i$ for every $1\le i \le r$, under the assumption
that $m_i \le \sum_{j\ne i} m_j$ for every $1\le i \le r$.
\end{itemize}

Now we are going to start with the recursion.

The resultant $\mathcal{R} = {\rm Res}(F_0, F_1^{(1)}, \dots,
F_{n_1}^{(1)}, \dots, F_1^{(r)}, \dots, F_{n_r}^{(r)})$ we want to
compute in the proof of Theorem \ref{eq-quasieq} involves a family
of multihomogeneous polynomials satisfying the conditions in (I)
with $m_i:= n_i$ and $m:= n$. Note that, for a generic game to
have a non-empty set of quasi-equilibria, the inequalities $m_i
\le\sum_{j\ne i } m_j $ for every $1\le i\le r$ hold.

Let $G_0^{(0)}$ and $G_k^{(i)}$ $(1\le i \le r,\, 1\le k \le m_i)$
be a family of polynomials satisfying the conditions in (I) and let
$\cI:= \{(0,0)\}\cup \{ (i,k) : 1\le i \le r,\ 1\le k\le m_i\}$.

First, note that ${\rm rank}(\cL_\cI)= m = \#\cI -1$. Let $J$ be a
proper subset of $\cI$. If there exist $(i, k)$, $(i', k')$ in $J$
with $i\ne i'$, then ${\rm rank}(\cL_J) = m \ge \#J$ and the same
holds if $(0,0) \in J$. On the other hand, if $J \subset \{(i,k):
1\le k \le m_i\}$ for a fixed $i\ne 0$, then ${\rm rank}(\cL_J) =
\sum_{j\ne i } m_j \ge m_i \ge \# J$. Therefore, the set of all
supports is the unique essential subset. So, the resultant is not
constant and the following identity holds:
$${\rm Res}(G_0^{(0)},(G^{(i)}_k)_{1\le i \le r; 1\le k \le m_i})
 = \prod_{\xi \in V} g_0^{(0)}(\xi)
\prod_{1\le j\le r}{\rm Res}(G^{(i)}_{kj})_{1\le i \le r, 1\le
k\le m_i}),$$ where $g_0^{(0)}$ is the dehomogeneized polynomial
obtained from $G_0^{(0)}$ by evaluating $x_{\ell\, m_\ell} = 1$
$(1 \le \ell \le r)$, $V$ is the set of common zeros in
$\mathbb{A}^m$ of the polynomials $g_k^{(i)}$ obtained in the same
way from the $G_k^{(i)}$ $(1\le i \le r, 1\le k \le m_i)$, and,
for each $1 \le j \le r$, $G^{(i)}_{kj}$ is the polynomial
obtained from $G^{(i)}_{k}$ by setting $x_{j\, m_j} = 0$.

\begin{itemize}
\item[(I.a)] Note that if $m_i \ge 2$ for every $1\le i \le r$,
(up to renaming variables and polynomials) each of the resultants
${\rm Res}(G^{(i)}_{kj})$ involves a family of polynomials
satisfying the conditions in (III).

\item[(I.b)] Without loss of generality, assume now that $m_1 =1$.
Here, we can discard the first group of variables. Then, we have
to compute the resultant of $m$ polynomials in $r-1$ groups of
$m_2+1, \dots, m_{r}+1$ variables each with $m_i$ polynomials with
multidegree $(1,\dots, 0, \dots, 1)$ where the zero is in the
$(i-1)$th coordinate for $2\le i \le r$ and one with multidegree
$(1,\dots, 1)$.

If, for every $2\le i \le r$, $m_i < \sum_{1\le j\le r,\, j\ne i}
m_j$, since $m_1 =1$ we deduce that $m_i \le \sum_{2\le j\le r,\,
j\ne i} m_j$ and so, the polynomial system obtained is of the form
(I) but with one group of variables less than the original one.

On the other hand, if $m_i = 1+\sum_{2\le j\le r,\, j\ne i} m_j$
for some $2\le i \le r$, we have that $m_j<m_i$ for every $j\ne
i$. Therefore, the unique essential subset is $\{(i, k): 1\le k\le
m_i\}$ and the resultant to be computed is the resultant of the
corresponding family of $m_i$ polynomials of multidegree
$(1,\dots, 1)$ in $r-2$ groups of $m_2+1, \dots, m_{i-1}+1,
m_{i+1}+1, \dots, m_{r}+1$ variables each, which is the situation
in (II).
\end{itemize}

If the polynomials involved satisfy the conditions in (II), all
the coordinates of their multidegrees are not zero and so, we can
apply the algorithm in \cite{JeSa} for the computation of their
resultant.

To analyze (III), let us consider first the case when $r=2$.
Here, the assumption on the numbers $m_i$ implies that $m_1 =
m_2:=M$.
\begin{enumerate}
\item[(III.a)] We consider the resultant of $M$ polynomials with
multidegrees $(0,1)$ and $M$ polynomials with multidegrees $(1,0)$
in two groups of $M$ and $M+1$ variables respectively. Now the
unique essential set is the corresponding to the first $M$
polynomials and, therefore their resultant equals the determinant
of the matrix of their coefficients.
\end{enumerate}

Assume now that $r>2$.
 Note that the equality $m_i = \sum_{j\ne i}
m_j$ may be valid for at most one value $i$. If, on the contrary,
$m_{i_1} = \sum_{j\ne {i_1}} m_j$ and $m_{i_2} = \sum_{j\ne i_2}
m_j$ hold for $i_1 \ne i_2$, it follows that $\sum_{j\ne i_1, i_2}
m_j =0$, which implies $r=2$. Let $G_k^{(i)}$ $(1\le i \le r,\,
1\le k \le m_i)$ be a family of polynomials satisfying the
conditions in (III).
\begin{enumerate}
\item[(III.b)] If $m_1=1$, we are under the same assumptions as in
(I.b). \item[(III.c)] If $m_1\ge 2$ and $m_i = \sum_{j\ne i} m_j$
for some $2\le i \le r$, the set $\{(i,k):1\le k \le m_i\}$ is
essential. The assumption $r>2$ implies that it is the unique
essential subset. Then, the resultant we compute is that of $m_i$
multilinear polynomials of multidegrees $(1,\dots, 1)$ in $r-1$
groups of $m_1, m_2+1, \dots, m_{i-1}+1, m_{i+1}+1\dots,m_r+1$
variables respectively. Thus, we are in situation (II).
\item[(III.d)] If $m_1 \ge 2$ and $m_i<\sum_{j\ne i} m_j$ for
every $2\le i \le r$, we deduce that $m_i\le m_1-1+\sum_{j\ne 1,i}
m_j$ for every $2\le i \le r$. Therefore, the unique essential
subset is the whole family of polynomials and applying Poisson's
formula we obtain:
$${\rm Res}((G^{(i)}_k)_{1\le i \le r; 1\le k \le m_i}) = \prod_{\xi \in V} g^{(1)}_{1}(\xi)
\prod_{2\le l\le r}{\rm Res}((G^{(1)}_{kl})_{2\le k\le
m_1};(G^{(i)}_{kl})_{2\le i \le r, 1\le k\le m_i}),$$ where
$g^{(1)}_{1}$ is the dehomogeneized of $G^{(1)}_{1}$ by setting
$x_{j m_j} = 1$ for every $1\le j \le r$, $V$ is the set of common
zeros in $\mathbb{A}^{m-1}$ of the polynomials $g^{(1)}_k$ $(2\le
k \le m_1)$, $g^{(i)}_k$ $(2\le i \le r, 1\le k \le m_i)$ obtained
in the same way from the homogeneous polynomias $G_k^{(i)}$,  and
$G^{(i)}_{kl}$ is the polynomial obtained from $G^{(i)}_{k}$ by
setting $x_{l\, m_l} = 0$.

According to the previous formula, we will have to compute $r-1$
multihomogeneous resultants of $m'=m-1$ polynomials each:

For $l=2,\dots, r$,  setting $m':= m-1$, $m'_1:= m_1-1$ and
$m'_i:= m_i$ for $i\ne 1$, the resultant to be computed involves
$m'$ polynomials in $r$ groups of $m'_1+1, \dots, m'_{l-1}+1,m'_l,
m'_{l+1}+1, \dots, m'_r+1$ variables each with $m'_i$ polynomials
of multidegree $d_i$ for every $1\le i \le r$. We have $m'_1 \le
\sum_{j\ne 1} m'_j$ and, for $i\ne 1$, the condition
$m_i<\sum_{j\ne i} m_j$ implies that $m'_i\le\sum_{j\ne i} m'_j$;
therefore, renaming variables and polynomials, we are again under
the assumptions of (III).

\end{enumerate}

\subsection{A bound for the degree of the resultant}

This section is devoted to proving an upper bound for the degree
$D$ of the resultant of multi-linear polynomials appearing in our
previous computations in terms of the total number $n$ of
strategies available to the players and the number $\delta$ of
totally mixed Nash equilibria of a generic game with the given
structure.

As before, we assume $n = n_1 +\cdots +n_r$ with $n_i \in \N$ for
every $1\le i \le r$. We consider the resultant $\cR$ of a family
of $n+1$ multi-linear polynomials in $r$ groups of $n_1+1, \dots,
n_r+1$ variables each, consisting of a polynomial $F_0$ of
multi-degree $d_0:= (1,\dots, 1)$ and, for every $1\le i \le r$, a
set of $n_i$ polynomials $F_k^{(i)}$ $(1\le k \le n_i)$ of
multi-degree $d_i:= (1,\dots, 1,0,1,\dots,1)$, where the $0$ lies
in the $i$th coordinate.

Throughout this section, we will denote
$$\delta_i:= {\rm
Bez}_{n_1,\dots, n_r}(d_0,1; d_{1}, n_{1}; \dots; d_{i},n_i-1;
\dots; d_{r}, n_r), \qquad i = 1,\dots, r,$$ the number of
solutions to a generic polynomial system of the indicated
multi-degrees. Recall that $\delta =  {\rm Bez}_{n_1,\dots, n_r}
(d_{1}, n_{1}; \dots; d_{r}, n_r)$.

The resultant $\cR$ is a multihomogeneous polynomial in the
coefficients of $F_0, \, F_k^{(i)}$ of degree $\delta$ in the
coefficients of $F_0$ and $\delta_i$ in the coefficients of
$F_k^{(i)}$ for every $1\le i \le r,\, 1\le k \le n_i$ (see, for
instance, \cite{PS}). Therefore, the total degree of $\cR$ equals
\begin{equation}\label{degR}
D = \delta + \sum_{1\le i\le r} n_i \delta_i.
\end{equation}
The following result provides an upper bound for $\delta_i$ which
will allow us to deduce an upper bound for $D$.

\begin{proposition}
For every $1\le i \le r$, we have $\delta_i\le (n_i +1) \delta$.
\end{proposition}

\begin{proof}{Proof.} Without loss of generality, we will prove
the stated upper bound for $i=1$.

First, note that the formula in the right hand side of identity
(\ref{bezoutmultihom}) for the computation of the B\' ezout number
of a multihomogeneous system of polynomial equations shows that
this number is additive in each of the multi-degrees involved.
Hence, as $\delta_1= {\rm Bez}_{n_1,\dots, n_r}(d_0,1; d_{1},
n_{1} -1; d_2, n_2;\dots;d_{r}, n_r)$ with $d_0 = d_1 +e_1$, where
$e_1 = (1,0,\dots, 0)$, we have
\begin{eqnarray*}
\delta_1 &=& {\rm Bez}_{n_1,\dots, n_r}(d_1, n_1;d_2, n_2; \dots;
d_r, n_r)+{\rm Bez}_{n_1,\dots, n_r}(e_1,1;d_1, n_1-1;d_2, n_2;
\dots; d_r, n_r)\\
&=& \delta +{\rm Bez}_{n_1,\dots, n_r}(e_1,1;d_1, n_1-1;d_2, n_2;
\dots; d_r, n_r).
\end{eqnarray*}
Now, identity (\ref{bezoutmultihom}) implies that $ {\rm
Bez}_{n_1,\dots, n_r}(e_1,1;d_1, n_1-1;d_2, n_2; \dots; d_r, n_r)
= \#\J_1,$ where
$$\J_1 = \{ (j_{11},\dots, j_{rn_r}) \, /\, j_{11} = 1,
j_{ik} \ne i \, \forall\, (i,k) \ne (1,1) \hbox{ and } \#\{ j_{hk}
\, /\, j_{hk}= i \} = n_i \ \forall \, 1\le i \le r\}.$$

In order to finish the proof, we will show that $\# \J_1 \le n_1
\delta$. Since $\delta$ equals the cardinality of the set $\J_0$
introduced in (\ref{J0}), we will compare the cardinalities of
both sets $\J_1$ and $\J_0$. To this end, we define the following
map from $\J_1$ to $\J_0$: with a given $n$-tuple $j:=(1, j_{12},
\dots, j_{1n_1}, \dots, j_{r1}, \dots, j_{rn_r}) \in \J_1$ we
associate the $n$-tuple $j' \in \J_0$ which is obtained by
exchanging the first coordinate of $j$ (which equals $1$) with the
first one which is different from $1$ and is located beyond
the $n_1$th coordinate. 
Note that a necessary condition for two distinct $n$-tuples in
$\J_1$ to lead to the same $n$-tuple in $\J_0$ by means of this
assignment is that they coincide in all of their coordinates
except for two of them located among the $n_1$ coordinates
$n_1+1,\dots, 2n_1$. Moreover,  the vector consisting of these
$n_1$ coordinates must be of the form $(1,\dots, 1, j_{hk},\dots)$
for both of them (possibly with no $1$ at the beginning) and so,
they can only differ in the length of the string of $1$'s in this
vector, which ranges between $0$ and $n_1 - 1$.  We conclude that
each element of $\J_0$ is the image of at most $n_1$ elements of(1) 
$\J_1$. It follows that $\# \J_1 \le n_1 \# \J_0$.

Therefore, we have that $\delta_1 \le \delta + n_1 \delta  =
(n_1+1) \delta $.
\end{proof}

Using the previous result along with identity (\ref{degR}) for the
degree $D$ of the resultant, we conclude:

\begin{corollary}\label{ubdegR}
With the previous assumptions and notations, $$D \le
\Big(1+\sum_{1\le i \le r} n_i (n_i+1) \Big) \delta \le n^2
\delta.$$
\end{corollary}

\begin{remark}
The bound stated in Corollary \ref{ubdegR} shows that all the
algorithms presented in this paper are polynomial in the number of
strategies $n_1, \dots, n_r$ of the $r$ players, and the generic
number $\delta$ of totally mixed Nash equilibria of a game with
the considered structure.
\end{remark}

\end{document}